\documentclass[aos,MSNbibl,seceqn,dvips]{arximspdf}
\usepackage{dcolumn}
\usepackage{graphicx}

\doi{10.1214/13-AOS1164} %kopijuoti is PTS
\volume{41}
\issue{5}
\pubyear{2013}
\firstpage{2608}
\lastpage{2638}

\makeatletter
\newcolumntype{d}[1]{D{.}{.}{#1}}

\newcommand{\rrvert}{\vert}
\newcommand{\llvert}{\vert}
\renewcommand{\Pr}{P}
\newtheorem{Theorem}{Theorem}[section]
\newtheorem{Lemma}[Theorem]{Lemma}
\newtheorem{Corollary}[Theorem]{Corollary}
\newtheorem{Proposition}[Theorem]{Proposition}
\newproclaim{Assumption}{Assumption}
\newproclaim{Remark}{Remark}[section]
\newproclaim{Example}{Example}[section]

\makeatother

\begin{document}
\begin{frontmatter}

\title{Local and global asymptotic inference in~smoothing spline models}
\runtitle{Asymptotic inference for smoothing spline}

\begin{aug}
\author[A]{\fnms{Zuofeng} \snm{Shang}\ead[label=e1]{zuofengshang@gmail.com}}
\and
\author[B]{\fnms{Guang} \snm{Cheng}\corref{}\ead[label=e2]{chengg@stat.purdue.edu}\thanksref{t2}}
\runauthor{Z. Shang and G. Cheng}
\affiliation{University of Notre Dame and Purdue University}
\address[A]{University of Notre Dame\\
Notre Dame, Indiana 46556\\
USA\\
\printead{e1}} %adresu isvedimo komanda gale!
\address[B]{Department of Statistics\\
Purdue University\\
250 N. University St.\\
West Lafayette, Indiana 47906\\
USA\\
\printead{e2}}
\end{aug}
\thankstext{t2}{Sponsored by NSF Grant DMS-09-06497 and CAREER Award DMS-11-51692.}

% HISTORY:
\received{\smonth{8} \syear{2013}}

% ABSTRACT
%
\begin{abstract}
This article studies local and global inference for smoothing spline
estimation in a unified asymptotic framework. We first introduce a new
technical tool called functional Bahadur representation, which
significantly \mbox{generalizes} the traditional Bahadur representation in
parametric models, that is, Bahadur [\textit{Ann. Inst. Statist. Math.}
\textbf{37} (1966) 577--580]. Equipped with this tool, we develop four
interconnected procedures for inference: (i)~pointwise confidence
interval; (ii)~local likelihood ratio testing; (iii)~simultaneous
confidence band; (iv)~global likelihood ratio testing. In particular,
our confidence intervals are proved to be asymptotically valid at any
point in the support, and they are shorter on average than the Bayesian
confidence intervals proposed by Wahba [\textit{J.~R. Stat. Soc. Ser. B Stat. Methodol.}
\textbf{45} (1983) 133--150] and Nychka [\textit{J.~Amer.
Statist. Assoc.} \textbf{83} (1988) 1134--1143]. We also discuss a
version of the Wilks phenomenon arising from local/global likelihood
ratio testing. It is also worth noting that our simultaneous confidence
bands are the first ones applicable to general quasi-likelihood models.
Furthermore, issues relating to optimality and efficiency are carefully
addressed. As a by-product, we discover a surprising relationship
between periodic and nonperiodic smoothing splines in terms of
inference.
\end{abstract}

% KEYWORDS
% Pirmas kwd is didziosios raides
%
\begin{keyword}[class=AMS]
\kwd[Primary ]{62G20}
\kwd{62F25}
\kwd[; secondary ]{62F15}
\kwd{62F12}
\end{keyword}
\begin{keyword}
\kwd{Asymptotic normality}
\kwd{functional Bahadur representation}
\kwd{local/global likelihood ratio test}
\kwd{simultaneous confidence band}
\kwd{smoothing spline}
\end{keyword}

\end{frontmatter}

%s1 #&#
\section{Introduction}\label{sec1}

As a flexible modeling tool, smoothing splines provide a general
framework for statistical analysis in a variety of fields; see
\cite{W90,W11,G02}. The asymptotic studies on smoothing splines in the
literature focus primarily on the estimation performance, and in
particular the global convergence. However, in practice it is often of
great interest to conduct \textit{asymptotic inference} on the unknown
functions. The procedures for inference developed in this article,
together with their rigorously derived asymptotic properties, fill this
long-standing gap in the smoothing spline literature.

As an illustration, consider two popular nonparametric regression
models: (i)~normal regression: $Y\mid Z=z\sim N(g_0(z),\sigma^2)$ for
some unknown $\sigma^2>0$; (ii)~logistic regression: $P(Y=1\mid
Z=z)=\exp (g_0(z))/(1+\exp(g_0(z)))$. The function $g_0$ is assumed to
be smooth in both models. Our goal in this paper is to develop
asymptotic theory for constructing pointwise confidence intervals and
simultaneous confidence bands for $g_0$, testing on the value of
$g_0(z_0)$ at any given point $z_0$, and testing whether $g_0$
satisfies certain global properties such as linearity. Pointwise
confidence intervals and tests on a local value are known as local
inference. Simultaneous confidence bands and tests on a global property
are known as global inference. To the best of our knowledge, there has
been little systematic and rigorous theoretical study of asymptotic
inference. This is partly because of the \mbox{technical} restrictions of the
widely used equivalent kernel method. The \textit{functional Bahadur
representation} (FBR) developed in this paper makes several important
contributions to this area. Our main contribution is a set of
procedures for local and global inference for a univariate smooth
function in a general class of nonparametric regression models that
cover both the aforementioned cases. Moreover, we investigate issues
relating to optimality and efficiency that have not been treated in the
existing smoothing spline literature.

The equivalent kernel has long been used as a standard tool for
handling the asymptotic properties of smoothing spline estimators, but
this method is restricted to least square regression; see
\cite{S84,MG93}. Furthermore, the equivalent kernel only
``approximates'' the reproducing kernel function, and the approximation
formula becomes extremely complicated when the penalty order increases
or the design points are nonuniform. To analyze the smoothing spline
estimate in a more effective way, we employ empirical process theory to
develop a new technical tool, the functional Bahadur representation,
which directly handles the ``exact'' reproducing kernel, and makes it
possible to study asymptotic inference in a broader range of
nonparametric models. An immediate consequence of the FBR is the
asymptotic normality of the smoothing spline estimate. This naturally
leads to the construction of pointwise asymptotic confidence intervals
(CIs). The classical Bayesian CIs in the literature \cite{W83,N88} are
valid on average over the observed covariates. However, our CIs are
proved to be theoretically valid at any point, and they even have
shorter lengths than the Bayesian CIs. We next introduce a likelihood
ratio method for testing the local value of a regression function. It
is shown that the null limiting distribution is a scaled Chi-square
with one degree of freedom, and that the scaling constant converges to
one as the smoothness level of the regression function increases.
Therefore, we have discovered an interesting Wilks phenomenon (meaning
that the asymptotic null distribution is free of nuisance parameters)
arising from the proposed nonparametric local testing.

Procedures for global inference are also useful. Simultaneous
confidence bands (SCBs) accurately depict the global behavior of the
regression function, and they have been extensively studied in the
literature. However, most of the efforts were devoted to simple
regression models with additive Gaussian errors based on kernel or
local polynomial estimates; see \cite{Har89,SL94,CK03,FZ00,ZP10}. By
incorporating the reproducing kernel Hilbert space (RKHS) theory into
\cite{BR73}, we obtain an SCB applicable to general nonparametric
regression models, and we demonstrate its theoretical validity based on
strong approximation techniques. To the best of our knowledge, this is
the first SCB ever developed for a general nonparametric regression
model in smoothing spline settings. We further demonstrate that our SCB
is optimal in the sense that the minimum width of the SCB achieves the
lower bound established by~\cite{GW08}. Model assessment is another
important aspect of global inference. Fan et al. \cite{FZZ01} used
local polynomial estimates for testing nonparametric regression models,
namely the generalized likelihood ratio test (GLRT). Based on smoothing
spline estimates, we propose an alternative method called the penalized
likelihood ratio test (PLRT), and we identify its null limiting
distribution as nearly Chi-square with diverging degrees of freedom.
Therefore, the Wilks phenomenon holds for the global test as well. More
importantly, we show that the PLRT achieves the minimax rate of testing
in the sense of \cite{I93}. In comparison, other smoothing-spline-based
tests such as the locally most powerful (LMP) test, the generalized
cross validation (GCV) test and the generalized maximum likelihood
ratio (GML) test (see \cite{CKWY88,W90,J96,C94,RG00,LW02}) either lead
to complicated null distributions with nuisance parameters or are not
known to be optimal.

As a by-product, we derive the asymptotic equivalence of the proposed
procedures based on periodic and nonperiodic smoothing splines under
mild conditions; see Remark~\ref{remark:equivker:appl}. In general,
our findings reveal an intrinsic connection between the two rather
different basis structures, which in turn facilitates the
implementation of inference.

Our paper is mainly devoted to theoretical studies. However, a few
practical issues are noteworthy. GCV is currently used for the
empirical tuning of the smoothing parameter, and it is known to result
in biased estimates if the true function is spatially inhomogeneous
with peaks and troughs. Moreover, the penalty order is prespecified
rather than data-driven. Future research is needed to develop an
efficient method for choosing a suitable smoothing parameter for bias
reduction and an empirical method for quantifying the penalty order
through data. We also note that some of our asymptotic procedures are
not fully automatic since certain quantities need to be estimated; see
Example~\ref{eglog}. A large sample size may be necessary for the
benefits of our asymptotic methods to become apparent. Finally, we want
to mention that extensions to more complicated models such as
multivariate smoothing spline models and semiparametric models are
conceptually feasible by applying similar FBR techniques and
likelihood-based approaches.

The rest of this paper is organized as follows.
Section~\ref{prelim:sec} introduces the basic notation, the model
assumptions, and some preliminary RKHS results. Section~\ref{secasy}
presents the FBR and the local asymptotic results. In
Sections~\ref{secjse} and~\ref{secjse2}, several procedures for local
and global inference together with their theoretical properties are
formally discussed. In Section~\ref{secexa}, we give three concrete
examples to illustrate our theory. Numerical studies are also provided
for both periodic and nonperiodic splines. The proofs are included in
an online supplementary document \cite{SC12}.

%s2 #&#
\section{Preliminaries}\label{prelim:sec}
%s2.1 #&#
\subsection{Notation and assumptions}\label{notation:assumption}

Suppose that the data $T_i=(Y_i,Z_i)$, $i=1,\ldots,n$, are i.i.d.
copies of $T=(Y,Z)$, where $Y\in\mathcal{Y}\subseteq\mathbb{R}$ is the
response variable, $Z\in\mathbb{I}$ is the
covariate variable and $\mathbb{I}=[0,1]$. Consider a general class of
nonparametric regression models under the primary assumption
%
%e2.1 #&#
\begin{equation}
\label{part:linear:model} \mu_0(Z)\equiv E(Y\mid Z)=F
\bigl(g_0(Z)\bigr),
\end{equation}
where $g_0(\cdot)$ is some unknown smooth function and $F(\cdot)$ is a
known link function. This framework covers two subclasses of
statistical interest. The first subclass assumes that the data are
modeled by $y_i\mid  z_i\sim p(y_i;\mu_0(z_i))$ for a conditional
distribution $p(y;\mu_0(z))$ unknown up to $\mu_0$. Instead of assuming
known distributions, the second subclass specifies the relation between
the conditional mean and variance as $\operatorname{Var}(Y\mid
Z)=\mathcal V(\mu_0(Z))$, where $\mathcal{V}$ is a known
positive-valued function. The nonparametric estimation of $g$ in the
second situation uses the quasi-likelihood $Q(y;
\mu)\equiv\int_{y}^\mu(y-s)/\mathcal V(s)\,ds$ as an objective function
(see \cite{W74}), where $\mu=F(g)$. Despite distinct modeling
principles, the two subclasses have a large overlap since $Q(y;\mu)$
coincides with $\log p(y;\mu)$ under many common combinations of $(F,
\mathcal V)$; see Table~2.1 of \cite{MN89}.

From now on, we focus on a smooth criterion function $\ell(y;a)\dvtx
\mathcal Y\times\mathbb R\mapsto\mathbb R$ that covers the above two
cases, that is, $\ell(y;a)=Q(y;F(a))$ or $\log p(y; F(a))$. Throughout
this paper, we define the functional parameter space $\mathcal H$ as
the \mbox{$m$th-}order Sobolev space:
\begin{eqnarray*}
H^m(\mathbb{I}) &\equiv& \bigl\{g \dvtx \mathbb{I} \mapsto
\mathbb{R}\mid  g^{(j)} \mbox{ is absolutely continuous }
\\
&&\hspace*{5pt}
\mbox{for }j=0,1,\ldots,m-1 \mbox{ and }g^{(m)} \in L_{2}(\mathbb{I})\bigr\},
\end{eqnarray*}
where $m$ is assumed to be known and larger than $1/2$. With some abuse
of notation, $\mathcal H$ may also refer to the homogeneous subspace
$H_0^m(\mathbb{I})$ of $H^m(\mathbb{I})$. The space $H_0^m(\mathbb
{I})$ is also known as the class of periodic functions such that a
function $g\in H_0^m(\mathbb{I})$ has the additional restrictions
$g^{(j)}(0)=g^{(j)}(1)$ for $j=0,1,\ldots,m-1$. Let $J(g,\widetilde
g)=\int_\mathbb{I} g^{(m)}(z)\widetilde g{}^{(m)}(z)\,dz$. Consider\vadjust{\goodbreak} the
penalized nonparametric estimate $\widehat g_{n,\lambda}$:
%
%e2.2 #&#
\begin{eqnarray}\label{est:f}
\widehat{g}_{n,\lambda}&=&\mathop{\arg\max}_{g\in\mathcal{H}}\ell
_{n,\lambda}(g)
\nonumber\\[-8pt]\\[-8pt]
&=& \mathop{\arg\max}_{g\in\mathcal{H}} \Biggl\{ \frac{1}{n}
\sum_{i=1}^n \ell\bigl(Y_i;g(Z_i)\bigr)-(\lambda/2) J(g,g) \Biggr\},\nonumber
\end{eqnarray}
where $J(g,g)$ is the roughness penalty and $\lambda$ is the smoothing
parameter, which converges to zero as $n\rightarrow\infty$. We use
$\lambda/2$ (rather than $\lambda$) to simplify future expressions. The
existence and uniqueness of $\widehat g_{n,\lambda}$ are guaranteed by
Theorem~2.9 of \cite{G02} when the null space $\mathcal
N_m\equiv\{g\in\mathcal{H}\dvtx J(g,g)=0\}$ is finite dimensional and
$\ell(y;a)$ is concave and continuous w.r.t. $a$.

We next assume a set of model conditions. Let $\mathcal{I}_0$ be the
range of $g_0$, which is obviously compact. Denote the first-, second-
and third-order derivatives of $\ell(y;a)$ w.r.t. $a$ by $\dot\ell
_a(y;a)$, $\ddot{\ell}_a(y;a)$ and $\ell'''_a(y;a)$, respectively.
We assume the following smoothness and tail conditions on $\ell$:

%
%asA.1 #&#
\begin{Assumption}\label{A.1}
\textup{(a)} $\ell(y;a)$ is three times continuously
    differentiable and concave w.r.t. $a$. There exists a bounded open
    interval $\mathcal{I}\supset\mathcal{I}_0$ and positive constants
    $C_0$ and $C_1$ s.t.
%
%e2.3 #&#
\begin{equation}
\label{A1:a:eq1} E \Bigl\{\exp\Bigl(\sup_{a\in\mathcal{I}}\bigl|\ddot{\ell
}_a(Y;a)\bigr|/C_0\Bigr) \bigm|Z \Bigr\}\leq C_0\qquad\mbox{a.s.}
\end{equation}
and
%
%e2.4 #&#
\begin{equation}
\label{A1:a:eq2} E \Bigl\{\exp\Bigl(\sup_{a\in\mathcal{I}} \bigl|
\ell'''_a(Y;a)\bigr|/C_0 \Bigr)\bigm| Z \Bigr\}\leq C_1\qquad\mbox{a.s.}
\end{equation}

\textup{(b)} There exists a positive constant $C_2$ such
    that $C_2^{-1}\leq I(Z)\equiv -E(\ddot{\ell}_a(Y;\break g_0(Z))\mid Z)\leq C_2$
    a.s.

\textup{(c)} $\epsilon\equiv \dot{\ell}_a(Y;g_0(Z))$ satisfies
$E(\epsilon\mid Z)=0$ and $E(\epsilon^2\mid Z)= I(Z)$ a.s.
\end{Assumption}

Assumption \ref{A.1}(a) implies the slow diverging rate $O_P(\log{n})$
of
\[
\max_{1\leq i\leq n}\sup_{a\in\mathcal{I}} \bigl|\ddot{\ell}_a(Y_i;a)\bigr|\quad\mbox{and}\quad
\max_{1\leq i\leq n}\sup_{a\in\mathcal {I}}\bigl|\ell'''_a(Y_i;a)\bigr|.
\]
When $\ell(y;a)=\log p(y;a)$, Assumption \ref{A.1}(b) imposes
boundedness and positive definiteness of the Fisher information, and
Assumption \ref{A.1}(c) trivially holds if $p$ satisfies certain
regularity conditions. When $\ell(y;a)=Q(y;\break F(a))$, we have
%
%e2.5 #&#
\begin{equation}\label{ddot:ell}
\ddot{\ell}_a(Y;a)=F_1(a)+\varepsilon F_2(a)\quad\mbox{and}\quad
\ell'''_a(Y;a)=\dot{F}_1(a)+\varepsilon\dot{F}_2(a),
\end{equation}
where\vspace*{1pt} $\varepsilon=Y-\mu_0(Z)$, $F_1(a)=-|\dot{F}(a)|^2/\mathcal
{V}(F(a))+(F(g_0(Z))-F(a))F_2(a)$ and
$F_2(a)=(\ddot{F}(a)\mathcal{V}(F(a))-\dot{\mathcal{V}}(F(a))|\dot
{F}(a)|^2)/\mathcal{V}^{2}(F(a))$. Hence, Assumption~\ref{A.1}(a) holds
if $F_j(a)$, $\dot{F}_j(a)$, $j=1,2$, are all bounded over
$a\in\mathcal{I}$ and
%
%e2.6 #&#
\begin{equation}
\label{a:condition:enno} E\bigl\{\exp\bigl(|\varepsilon|/C_0\bigr) \mid  Z\bigr\}
\leq C_1\qquad\mbox{a.s.}
\end{equation}
By (\ref{ddot:ell}), we have $I(Z)=|\dot{F}(g_0(Z))|^2/\mathcal
{V}(F(g_0(Z)))$. Thus, Assumption \ref{A.1}(b) holds if
%
%e2.7 #&#
\begin{equation}\label{quasi:Iu}
1/C_2\leq\frac{|\dot{F}(a)|^2}{\mathcal{V}(F(a))}\leq C_2\qquad\mbox{for all $a\in\mathcal{I}_0$ a.s.}
\end{equation}
Assumption \ref{A.1}(c) follows from the definition of $\mathcal
V(\cdot)$. The sub-exponential tail condition (\ref{a:condition:enno})
and the boundedness condition (\ref{quasi:Iu}) are very mild
quasi-likelihood model assumptions (e.g., also assumed in
\cite{MVG97}). The assumption that $F_j$ and $\dot{F}_j$ are both
bounded over $\mathcal{I}$ could be restrictive and can be removed in
many cases, such as the binary logistic regression model, by applying
empirical process arguments similar to those in Section~7 of
\cite{MVG97}.

%s2.2 #&#
\subsection{Reproducing kernel Hilbert space}\label{secRKHS}
We now introduce a number of RKHS results as extensions of \cite{CO90}
and \cite{N95}. It is well known that, when $m>1/2$, $\mathcal
{H}=H^m(\mathbb{I})$ [or $H_0^m(\mathbb{I})$] is an RKHS endowed with
the inner product $\langle g,\widetilde{g}\rangle=
E\{I(Z)g(Z)\widetilde{g}(Z)\}+\lambda J(g,\widetilde{g})$ and the norm
%
%e2.8 #&#
\begin{equation}
\label{norm:Hm} \|g\|^2=\langle g,g\rangle.
\end{equation}
The reproducing kernel $K(z_1,z_2)$ defined on
$\mathbb{I}\times\mathbb{I}$ is known to have the following property:
\[
K_z(\cdot)\equiv K(z,\cdot)\in\mathcal{H}\quad\mbox{and}\quad
\langle K_z,g\rangle=g(z)
\qquad\mbox{for any $z\in\mathbb{I}$ and $g\in \mathcal{H}$}.
\]
Obviously, $K$ is symmetric with $K(z_1,z_2)=K(z_2,z_1)$. We further
introduce a~positive definite self-adjoint operator $W_\lambda\dvtx
\mathcal{H}\mapsto\mathcal{H}$ such that
%
%e2.9 #&#
\begin{equation}
\label{sec2:eq:0} \langle W_\lambda g,\widetilde{g}\rangle=\lambda J(g,
\widetilde{g})
\end{equation}
for any $g,\widetilde{g}\in\mathcal{H}$. Let
$V(g,\widetilde{g})=E\{I(Z)g(Z)\widetilde{g}(Z)\}$. Then $\langle
g,\widetilde{g}\rangle=V(g,\widetilde{g})+\langle
W_\lambda g,\widetilde{g}\rangle$ and
$V(g,\widetilde{g})=\langle
(\mathrm{id}-W_\lambda)g,\widetilde{g}\rangle$,
where $\mathrm{id}$ denotes the identity operator.

Next, we assume that there exists a sequence of basis functions in the
space $\mathcal{H}$ that simultaneously diagonalizes the bilinear forms
$V$ and $J$. Such eigenvalue/eigenfunction assumptions are typical in
the smoothing spline literature, and they are critical to control the
local behavior of the penalized estimates. Hereafter, we denote
positive sequences $a_\mu$ and $b_\mu$ as $a_\mu\asymp b_\mu$ if they
satisfy $\lim_{\mu\rightarrow\infty}(a_\mu/b_\mu)=c>0$. If $c=1$, we
write $a_\mu\sim b_\mu$. Let $\sum_\nu$ denote the sum over
$\nu\in\mathbb{N}=\{0,1,2,\ldots\}$ for convenience. Denote the
sup-norm of $g\in\mathcal{H}$ as
$\|g\|_{\sup}=\sup_{z\in\mathbb{I}}|g(z)|$.

%
%asA.2 #&#
\begin{Assumption} \label{A.3}
There exists a sequence of eigenfunctions $h_\nu\in \mathcal{H}$
satisfying $\sup_{\nu\in\mathbb{N}}\|h_\nu\|_{\sup}<\infty$, and a
nondecreasing sequence of eigenvalues $\gamma_\nu\asymp\nu^{2m}$ such
that
%
%e2.10 #&#
\begin{equation}
\label{regular} V(h_\mu,h_\nu)=\delta_{\mu\nu},\qquad
J(h_\mu,h_\nu)=\gamma_\mu\delta_{\mu\nu},\qquad
\mu,\nu\in\mathbb{N},
\end{equation}
where $\delta_{\mu\nu}$ is the Kronecker's delta. In particular, any
$g\in\mathcal{H}$ admits a Fourier expansion
$g=\sum_{\nu} V(g,h_\nu)h_\nu$ with convergence in the $\|\cdot\|$-norm.
\end{Assumption}

Assumption~\ref{A.3} enables us to derive explicit expressions for $\|
g\|$, $K_z(\cdot)$ and $W_\lambda h_\nu(\cdot)$ for any $g\in
\mathcal{H}$ and $z\in\mathbb I$; see Proposition~\ref{Prop:RPK} below.

%
%pr2.1 #&#
\begin{Proposition}\label{Prop:RPK}
For any $g\in\mathcal{H}$ and $z\in\mathbb{I}$, we have $\|g\|
^2=\sum_{\nu}
|V(g,\break h_\nu)|^2(1+\lambda\gamma_\nu)$, $K_z(\cdot)=\sum_\nu
\frac{h_\nu(z)}{1+\lambda\gamma_\nu} h_\nu(\cdot)$ and
$W_\lambda h_\nu(\cdot)=
\frac{\lambda\gamma_\nu}{1+\lambda\gamma_\nu}h_\nu(\cdot)$
under Assumption~\ref{A.3}.
\end{Proposition}

For future theoretical derivations, we need to figure out the
underlying eigensystem that implies Assumption~\ref{A.3}. For example,
when $\ell(y;a)=-(y-a)^2/2$ and $\mathcal H=H_0^m(\mathbb I)$,
Assumption~\ref{A.3} is known to be satisfied if $(\gamma_\nu, h_\nu)$
is chosen as the trigonometric polynomial basis specified in
(\ref{H0:basis}) of Example~\ref{exampss}. For general $\ell(y;a)$ with
$\mathcal H=H^m(\mathbb{I})$, Proposition \ref{BVP} below says that
Assumption \ref{A.3} is still valid if $(\gamma_\nu, h_\nu)$ is chosen
as the (normalized) solution of the following equations:
%
%e2.11 #&#
%e2.12 #&#
\begin{eqnarray}\label{eigen:problem}
(-1)^m h_\nu^{(2m)}(\cdot)=\gamma_\nu I(\cdot)\pi(\cdot) h_\nu(\cdot),\qquad
h_\nu^{(j)}(0)=h_\nu^{(j)}(1)=0,
\nonumber\\[-8pt]\\[-8pt]
\eqntext{j=m,m+1,\ldots,2m-1,}
\end{eqnarray}
where $\pi(\cdot)$ is the marginal density of the covariate $Z$.
Proposition~\ref{BVP} can be viewed as a nontrivial extension of
\cite{U88}, which assumes $I=\pi=1$. The proof relies substantially on
the ODE techniques developed in \cite{Birk1908,Stone1926}. Let
$C^{m}(\mathbb I)$ be the class of the $m$th-order continuously
differentiable functions over $\mathbb I$.

%
%pr2.2 #&#
\begin{Proposition}\label{BVP}
If $\pi(z), I(z)\in C^{2m-1}(\mathbb{I})$ are both bounded away from
zero and infinity over $\mathbb{I}$, then the eigenvalues $\gamma_\nu$
and the corresponding eigenfunctions $h_\nu$, found from
(\ref{eigen:problem}) and normalized to $V(h_\nu,h_\nu)=1$, satisfy
Assumption~\ref{A.3}.
\end{Proposition}

Finally, for later use we summarize the notation for Fr\'{e}chet derivatives.
Let $\Delta g$, $\Delta g_{j}\in\mathcal{H}$ for $j=1,2,3$.
The Fr\'{e}chet derivative of $\ell_{n,\lambda}$ can be identified as
\begin{eqnarray*}
D\ell_{n,\lambda}(g)\Delta g&=&\frac{1}{n}\sum
_{i=1}^{n}\dot{\ell}_a
\bigl(Y_i; g(Z_i)\bigr) \langle K_{Z_i},\Delta
g\rangle-\langle W_\lambda g,\Delta g\rangle
\\
&\equiv&\bigl\langle S_n(g), \Delta g\bigr\rangle-\langle
W_\lambda g,\Delta g\rangle
\\
&\equiv&\bigl\langle S_{n,\lambda}(g),\Delta
g \bigr\rangle.
\end{eqnarray*}
Note that $S_{n,\lambda}(\widehat{g}_{n,\lambda})=0$ and
$S_{n,\lambda}(g_0)$ can be expressed as
%
%e2.13 #&#
\begin{equation}\label{score}
S_{n,\lambda}(g_0)=\frac{1}{n} \sum
_{i=1}^n\epsilon_i K_{Z_i}-W_\lambda g_0.
\end{equation}
The Fr\'{e}chet derivative of $S_{n,\lambda}$ $(DS_{n,\lambda})$ is
denoted $DS_{n,\lambda}(g)\Delta g_1\Delta g_2
(D^2S_{n,\lambda}(g)\*\Delta g_1\Delta g_2 \Delta g_3)$. These\vadjust{\goodbreak}
derivatives can be explicitly written as $D^2 \ell_{n,\lambda
}(g)\Delta g_1\Delta g_2=n^{-1}\sum_{i=1}^n
\ddot{\ell}_a(Y_i;g(Z_i))\langle K_{Z_i}, \Delta g_1\rangle\langle
K_{Z_i},\Delta g_2\rangle-\langle W_\lambda\Delta g_1,\Delta
g_2\rangle$ [or $D^3 \ell_{n,\lambda}(g)\*\Delta g_1\Delta g_2
\Delta g_3=n^{-1}\sum_{i=1}^n \ell'''_a(Y_i;g(Z_i))\langle K_{Z_i},
\Delta g_1\rangle\langle K_{Z_i},\Delta g_2\rangle\langle
K_{Z_i},\Delta g_3\rangle$].

Define $S(g)=E\{S_n(g)\}$, $S_\lambda(g)=S(g)-W_\lambda g$ and
$DS_\lambda(g)=DS(g)-W_\lambda$, where $DS(g)\Delta g_1\Delta
g_2=E\{\ddot{\ell}_a(Y; g(Z))\langle K_{Z}, \Delta g_1\rangle\langle
K_{Z},\Delta g_2\rangle\}$. Since $\langle DS_\lambda(g_0) f,\break
g\rangle=-\langle f,g\rangle$ for any $f,g\in\mathcal{H}$, we
have the following result:

%
%pr2.3 #&#
\begin{Proposition}\label{basic:prop}
$DS_\lambda(g_0)=-\mathrm{id}$, where $\mathrm{id}$ is the identity
operator on~$\mathcal{H}$.
\end{Proposition}

%s3 #&#
\section{Functional Bahadur representation}\label{secasy}
In this section, we first develop the key technical tool of this paper:
\textit{functional Bahadur representation}, and we then present the
local asymptotics of the smoothing spline estimate as a straightforward
application. In fact, FBR provides a rigorous theoretical foundation
for the procedures for inference developed in Sections~\ref{secjse}
and~\ref{secjse2}.

%s3.1 #&#
\subsection{Functional Bahadur representation}\label{secfbr}
We first present a relationship between the norms $\|\cdot\|_{\sup}$
and $\|\cdot\|$ in Lemma~\ref{lemma0} below, and we then derive a
\textit{concentration inequality} in Lemma~\ref{lemma2} as the
preliminary step for obtaining the FBR. For convenience, we denote
$\lambda^{1/(2m)}$ as $h$.

%
%le3.1 #&#
\begin{Lemma}\label{lemma0}
There exists a constant $c_m>0$ s.t. $|g(z)|\leq c_m h^{-1/2}\|g\|$ for
any $z\in\mathbb{I}$ and $g\in\mathcal{H}$. In particular, $c_m$ is not
dependent on the choice of $z$ and $g$. Hence, $\|g\|_{\sup}\leq c_m
h^{-1/2}\|g\|$.
\end{Lemma}

Define $\mathcal{G}=\{g\in\mathcal{H}\dvtx \|g\|_{\sup}\leq1,
J(g,g)\leq c_m^{-2} h\lambda^{-1}\}$, where $c_m$ is specified in Lemma
\ref{lemma0}. Recall that $\mathcal{T}$ denotes the domain of the full
data variable $T=(Y,Z)$. We now prove a concentration inequality on the
empirical process $Z_n(g)$ defined, for any $g\in\mathcal{G}$ and
$z\in\mathbb I$ as
%
%e3.1 #&#
\begin{equation}\label{Znf}
Z_n(g)
(z)=\frac{1}{\sqrt{n}}\sum_{i=1}^n \bigl[\psi_n(T_i;
g)K_{Z_i}(z)-E\bigl( \psi_n(T; g)K_Z(z)\bigr)\bigr],
\end{equation}
where $\psi_n(T; g)$ is a real-valued function (possibly
depending on $n$) defined on
$\mathcal{T}\times\mathcal{G}$.

%
%le3.2 #&#
\begin{Lemma}\label{lemma2}
Suppose that $\psi_n$ satisfies the following Lipschitz continuity condition:
%
%e3.2 #&#
\begin{equation}\label{Lip:cont:psi}
\qquad\bigl|\psi_n(T; f)-\psi_n(T; g)\bigr|\leq
c_m^{-1} h^{1/2} \|f-g\|_{\sup}\qquad \mbox{for any } f,g\in\mathcal{G},
\end{equation}
where $c_m$ is specified in Lemma \ref{lemma0}. Then we have
\[
\lim_{n\rightarrow\infty}P \biggl(\sup_{g\in\mathcal{G}}
\frac{\|
Z_n(g)\|}{h^{-(2m-1)/(4m)}\|g\|_{\sup}^{1-1/(2m)}+n^{-1/2}}\leq(5\log
\log{n})^{1/2} \biggr)=1.
\]
\end{Lemma}

To obtain the FBR, we need to further assume a proper convergence rate
for~$\widehat g_{n,\lambda}$:

%
%asA.3 #&#
\begin{Assumption}\label{A.5}
$\|\widehat{g}_{n,\lambda}-g_0\|=O_P((nh)^{-1/2}+h^m)$.
\end{Assumption}

Simple (but not necessarily the weakest) sufficient conditions for
Assumption~\ref{A.5} are provided in
Proposition~\ref{rates:convergence} below. Before stating this result,
we introduce another norm on the space $\mathcal{H}$, that is, more
commonly used in functional analysis. For any $g\in\mathcal{H}$, define
%
%e3.3 #&#
\begin{equation}
\label{another:sobnorm} \|g\|_{\mathcal{H}}^2=E\bigl
\{I(Z)g(Z)^2\bigr\}+J(g,g).
\end{equation}
When $\lambda\leq1$, $\|\cdot\|_{\mathcal H}$ is a type of Sobolev norm
dominating $\|\cdot\|$ defined in (\ref{norm:Hm}). Denote
%
%e3.4 #&#
\begin{equation}
\label{optsm} \lambda^\ast\asymp
n^{-2m/(2m+1)}\qquad \mbox{or equivalently, } h^\ast\asymp
n^{-1/(2m+1)}.
\end{equation}
Note that $\lambda^\ast$ is known
as the optimal order when we estimate $g_0\in\mathcal{H}$.

%
%pr3.3 #&#
\begin{Proposition}\label{rates:convergence}$\!\!\!\!\!$
Suppose that Assumption \ref{A.1} holds, and further $\|\widehat
{g}_{n,\lambda}-g_0\|_{\mathcal{H}}=o_P(1)$. If
$h$ satisfies $(n^{1/2}h)^{-1}(\log\log{n})^{m/(2m-1)}(\log
{n})^{2m/(2m-1)}=o(1)$, then Assumption~\ref{A.5} is satisfied. In
particular, $\widehat g_{n,\lambda}$ achieves the optimal rate of
convergence, that is, $O_P(n^{-m/(2m+1)})$, when
$\lambda=\lambda^\ast$.
\end{Proposition}

Classical results on rates of convergence are obtained through either
linearization techniques, for example, \cite{CO90}, or quadratic
approximation devices, for example, \cite{GQ93,G02}.
However, the proof of Proposition \ref{rates:convergence} relies on
empirical process techniques. Hence, it is not surprising that
Proposition \ref{rates:convergence} requires a different set of
conditions than those used in \cite{CO90,GQ93,G02}, although the
derived convergence rates are the same and in all approaches the
optimal rate is achieved when $\lambda=\lambda^\ast$. For example,
Cox and O'Sullivan \cite{CO90} assumed a weaker smoothness condition
on the likelihood function but a more restrictive condition on $h$,
that is, $(n^{1/2}h\lambda^\alpha)^{-1}=o(1)$ for some $\alpha>0$.

Now we are ready to present the key technical tool: \textit{functional
Bahadur representation}, which is also of independent interest. Shang
\cite{S10} developed a different form of Bahadur representation, which
is of limited use in practice. This is due to the intractable form of
the inverse operator $DS_\lambda(g_0)^{-1}$, constructed based on a
different type of Sobolev norm. However, by incorporating $\lambda$
into the norm (\ref{norm:Hm}), we can show $DS_\lambda(g_0)^{-1}=-\mathrm{id}$
based on Proposition \ref{basic:prop}, and thus obtain a more refined
version of the representation of \cite{S10} that naturally applies to
our general setting for inference purposes.

%
%th3.4 #&#
\begin{Theorem}[(Functional Bahadur representation)]\label{main:thm1}
Suppose that Assumptions \ref{A.1}--\ref{A.5} hold, $h=o(1)$ and
$nh^2\rightarrow\infty$. Recall that $S_{n,\lambda}(g_0)$ is defined
in~(\ref{score}). Then we have
%
%e3.5 #&#
\begin{equation}
\label{baha:reps} \bigl\|\widehat{g}_{n,\lambda}-g_0-S_{n,\lambda}(g_0)
\bigr\|=O_P(a_n\log{n}),
\end{equation}
where $a_n=n^{-1/2}((nh)^{-1/2}\hspace*{-0.2pt}+h^m)
h^{-(6m-1)/(4m)}(\log\log{n})^{1/2}\hspace*{-0.2pt}+C_\ell h^{-1/2}
((nh)^{-1}\hspace*{-0.6pt}+h^{2m})/\log{n}$ and $C_\ell=\sup_{z\in\mathbb{I}}E\{
\sup_{a\in\mathcal{I}}|\ell'''_a(Y;a)| \mid Z=z\}$. When $h=h^\ast
$, the RHS of (\ref{baha:reps}) is $o_{P}(n^{-m/(2m+1)})$.
\end{Theorem}

%s3.2 #&#
\subsection{Local asymptotic behavior}\label{sec3.2}
In this section, we obtain the pointwise asymptotics of $\widehat
g_{n,\lambda}$ as a direct application of the FBR. The equivalent
kernel method may be used for this purpose, but it is restricted to
$L_2$ regression, for example, \cite{S84}. However, the FBR-based
proof applies to more general regression. Notably, our results reveal
that several well-known global convergence properties continue to hold locally.

%
%th3.5 #&#
\begin{Theorem}[(General regression)]\label{main:thm2}
Assume Assumptions \ref{A.1}--\ref{A.5}, and suppose $h=o(1)$,
$nh^2\rightarrow\infty$ and $a_n\log{n}=o(n^{-1/2})$, where $a_n$ is
defined in Theorem~\ref{main:thm1}, as $n\rightarrow \infty$.
Furthermore, for any $z_0\in\mathbb{I}$,
%
%e3.6 #&#
\begin{equation}\label{main:thm2:part2}
h V(K_{z_0},K_{z_0})\rightarrow
\sigma^2_{z_0}\qquad\mbox{as } n\rightarrow\infty.
\end{equation}
Let $g_0^*=(\mathrm{id}-W_\lambda)g_0$ be the biased ``true parameter.'' Then we
have
%
%e3.7 #&#
\begin{equation}
\label{main:thm2:eq} \sqrt{nh}\bigl(\widehat{g}_{n,\lambda}(z_0)-g_0^*(z_0)
\bigr) \stackrel{d} {\longrightarrow} N\bigl(0,\sigma_{z_0}^2
\bigr),
\end{equation}
where
%
%e3.8 #&#
\begin{equation}
\label{weighted:sum:sigma} \sigma_{z_0}^2=\lim
_{h\rightarrow0}\sum_{\nu}\frac{h|h_\nu
(z_0)|^2}{(1+\lambda\gamma_\nu)^2}.
\end{equation}
\end{Theorem}

From Theorem \ref{main:thm2}, we immediately obtain the following result.

%
%co3.6 #&#
\begin{Corollary}\label{main:thm3}
Suppose that the conditions in Theorem \ref{main:thm2} hold and
%
%e3.9 #&#
\begin{equation}\label{weta}
\lim_{n\rightarrow\infty}(nh)^{1/2}(W_\lambda
g_0) (z_0)=-b_{z_0}.
\end{equation}
Then we have
%
%e3.10 #&#
\begin{equation}\label{an:important:limit}
\sqrt{nh}\bigl(\widehat{g}_{n,\lambda}(z_0)-g_0(z_0)
\bigr)\stackrel{d} {\longrightarrow} N\bigl(b_{z_0},
\sigma_{z_0}^2\bigr),
\end{equation}
where $\sigma_{z_0}^2$ is defined as in (\ref{weighted:sum:sigma}).
\end{Corollary}

To illustrate Corollary \ref{main:thm3} in detail, we consider $L_2$
regression in which $W_\lambda g_0(z_0)$ (also $b_{z_0}$) has an
explicit expression under the additional boundary conditions:
%
%e3.11 #&#
\begin{equation}
\label{BC:eta0} g^{(j)}_0(0)=g_0^{(j)}(1)=0\qquad\mbox{for } j=m,\ldots,2m-1.
\end{equation}
In fact, (\ref{BC:eta0}) is also the price we pay for obtaining the
boundary results, that is, $z_0=0, 1$. However, (\ref{BC:eta0}) could
be too strong in practice. Therefore, we provide an alternative set of
conditions in (\ref{limit:for:K}) below, which can be implied by the
so-called ``exponential envelope condition'' introduced in \cite{N95}.
In Corollary~\ref{l2:regression:CLT} below, we consider two different
cases: $b_{z_0}\neq0$ and $b_{z_{0}}=0$.

%
%co3.7 #&#
\begin{Corollary}[($L_2$ regression)]\label{l2:regression:CLT}
Let $m>(3+\sqrt{5})/4\approx1.309$ and\break $\ell(y;a)=-(y-a)^2/2$. Suppose
that Assumption \ref{A.5} and (\ref{main:thm2:part2}) hold, and the
normalized eigenfunctions $h_\nu $ satisfy (\ref{eigen:problem}).
Assume that $g_0\in H^{2m}(\mathbb{I})$ satisfies
$\sum_\nu|V(g_0^{(2m)},h_\nu)h_\nu(z_0)|<\infty$.
\begin{longlist}[(iii)]
\item[(i)] Suppose $g_0$ satisfies the boundary conditions
    (\ref{BC:eta0}).
    If $h/n^{-1/(4m+1)}\rightarrow c>0$, then we have, for any
    $z_0\in[0,1]$,
%
%e3.12 #&#
\begin{equation}\label{Cor:limit}
\qquad\sqrt{nh}
\bigl(\widehat{g}_{n,\lambda}(z_0)-g_0(z_0)
\bigr)\stackrel{d} {\longrightarrow} N \bigl((-1)^{m-1}c^{2m}g_0^{(2m)}(z_0)/
\pi(z_0),\sigma_{z_0}^2 \bigr).
\end{equation}
If $h\asymp n^{-d}$ for some $\frac{1}{4m+1}<d\leq
\frac{2m}{8m-1}$, then we have, for any $z_0\in[0,1]$,
%
%e3.13 #&#
\begin{equation}
\label{Cor:limit1}
\sqrt{nh}\bigl(
\widehat{g}_{n,\lambda}(z_0)-g_0(z_0)
\bigr)\stackrel{d} {\longrightarrow}N\bigl(0,\sigma_{z_0}^2
\bigr).
\end{equation}

\item[(ii)] If we replace the boundary conditions (\ref{BC:eta0}) by
    the following reproducing kernel conditions: for any $z_0\in(0,1)$, as
    $h\rightarrow0$
%
%e3.14 #&#
%e3.15 #&#
\begin{eqnarray}\label{limit:for:K}
\frac{\partial^{j}}{\partial z^j} K_{z_0}(z) \biggm|_{z=0}=o(1),\qquad
\frac{\partial^{j}}{\partial z^j} K_{z_0}(z) \biggm|_{z=1}=o(1)
\nonumber\\[-10pt]\\[-6pt]
\eqntext{\mbox{for } j=0,\ldots,m-1,}
\end{eqnarray}
then (\ref{Cor:limit}) and (\ref{Cor:limit1}) hold for any
$z_0\in(0,1)$.
\end{longlist}
\end{Corollary}

We note that in (\ref{Cor:limit}) the asymptotic bias is proportional
to $g_0^{(2m)}(z_0)$, and the asymptotic variance can be expressed as a
weighted sum of squares of the (\textit{infinitely many}) terms $h_\nu
(z_0)$; see (\ref{weighted:sum:sigma}). These observations are
consistent with those in the polynomial spline setting insofar as the
bias is proportional to $g_0^{(2m)}(z_0)$, and the variance is a
weighted sum of squares of (\,\textit{finitely many}) terms involving the
normalized B-spline basis functions evaluated at $z_0$; see
\cite{ZSW98}. Furthermore, (\ref{Cor:limit1})~describes how to remove
the asymptotic bias through \mbox{undersmoothing,} although the corresponding
smoothing parameter yields suboptimal estimates in terms of the
convergence rate.

The existing smoothing spline literature is mostly concerned with the
global convergence properties of the estimates. For example, Nychka
\cite{N95} and Rice and Rosenblatt \cite{RR83} derived global
convergence rates in terms of the (integrated) mean squared error.
Instead, Theorem \ref{main:thm2} and Corollaries \ref{main:thm3} and
\ref{l2:regression:CLT} mainly focus on local asymptotics, and they
conclude that the well-known global results on the rates of convergence
also hold in the \textit{local} sense.

%s4 #&#
\section{Local asymptotic inference}\label{secjse}
We consider inferring $g(\cdot)$ \textit{locally} by constructing the
pointwise asymptotic CI in Section~\ref{secconf} and testing the local
hypothesis in Section~\ref{seclrt}.

%s4.1 #&#
\subsection{Pointwise confidence interval}\label{secconf}
We consider the confidence interval for some real-valued smooth
function of $g_0(z_0)$ at any fixed $z_0\in\mathbb I$, denoted
$\rho_0=\rho(g_0(z_0))$, for example, $\rho_0=\exp
(g_0(z_0))/(1+\exp(g_0(z_0)))$ in logistic regression. Corollary
\ref{main:thm3} together with the Delta method immediately implies
Proposition~\ref{Prop:conf} on the pointwise CI where the asymptotic
estimation bias is assumed to be removed by undersmoothing.

%
%pr4.1 #&#
\begin{Proposition}[(Pointwise confidence interval)]\label{Prop:conf}
Suppose that the assumptions in Corollary
\ref{main:thm3} hold and that the estimation bias asymptotically
vanishes, that is, $\lim_{n\rightarrow\infty}(nh)^{1/2}(W_\lambda
g_0)(z_0)=0$. Let $\dot\rho(\cdot)$ be the first derivative of
$\rho(\cdot)$. If $\dot\rho(g_0(z_0))\neq0$, we have
\[
P \biggl(\rho_0\in\biggl[\rho\bigl(\widehat{g}_{n,\lambda}(z_0)
\bigr) \pm\Phi(\alpha/2)\frac{\dot\rho(g_0(z_0))\sigma_{z_0}}{\sqrt
{nh}} \biggr] \biggr)\longrightarrow1-
\alpha,
\]
where $\Phi(\alpha)$ is the lower $\alpha$th quantile of $N(0,1)$.
\end{Proposition}

From now on, we focus on the pointwise CI for $g_0(z_0)$ and compare it
with the classical \textit{Bayesian confidence intervals} proposed by
Wahba \cite{W83} and Nychka~\cite{N88}. For simplicity, we consider
$\ell(y;a)=-(y-a)^2/(2\sigma^2)$, $Z\sim \operatorname{Unif}[0,1]$ and
\mbox{$\mathcal H=H_0^{m}(\mathbb I)$} under which Proposition~\ref{Prop:conf}
implies the following asymptotic 95\% CI for $g_0(z_0)$:
%
%e4.1 #&#
\begin{equation}
\label{confintformula:refine} \widehat{g}_{n,\lambda}(z_0)
\pm1.96\sigma\sqrt{I_2/\bigl(n\pi h^\dag\bigr)},
\end{equation}
where $h^\dag=h\sigma^{1/m}$ and $I_l=\int_0^1 (1+x^{2m})^{-l}\,dx$ for
$l=1,2$; see case~(I) of Example~\ref{exampss} for the derivations.
When $\sigma$ is unknown, we may replace it by any consistent estimate.
As far as we are aware, (\ref{confintformula:refine}) is the first
rigorously proven pointwise CI for smoothing spline. It is well known
that the Bayesian type CI only approximately achieves the 95\% nominal
level on average rather than pointwise. Specifically, its average
coverage probability over the observed covariates is \textit{not}
exactly 95\% even asymptotically. Furthermore, the Bayesian type CI
ignores the important issue of coverage uniformity across the design
space, and thus it may not be reliable if only evaluated at peaks or
troughs, as pointed out in \cite{N88}. However, the asymptotic CI
(\ref{confintformula:refine}) is proved to be valid at any point, and
is even shorter than the Bayesian CIs proposed in \cite{W83,N88}.

As an illustration, we perform a detailed comparison of the three CIs
for the special case $m=2$. We first derive the asymptotically
equivalent versions of the Bayesian CIs. Wahba \cite{W83} proposed the
following heuristic CI
under a Bayesian framework:
%
%e4.2 #&#
\begin{equation}
\label{wahba:bci} \widehat{g}_{n,\lambda}(z_0)\pm1.96\sigma\sqrt{a
\bigl(h^\dag\bigr)},
\end{equation}
where $a(h^\dag)=n^{-1} (1+(1+(\pi n h^\dag))^{-4}+2\sum_{\nu
=1}^{n/2-1}(1+(2\pi\nu h^\dag))^{-4} )$. Under the assumptions
$h^\dag=o(1)$ and $(nh^\dag)^{-1}=o(1)$, Lemma~\ref{lemeg1} in Example
\ref{exampss} implies $2\sum_{\nu =1}^{n/2-1}(1+(2\pi\nu
h^\dag))^{-4}\sim I_1/(\pi h^\dag)=4I_2/(3\pi h^\dag)$, since
$I_2/I_1=3/4$ when $m=2$. Hence, we obtain an asymptotically equivalent
version of Wahba's Bayesian CI as
%
%e4.3 #&#
\begin{equation}
\label{asympwahba:bci} \widehat{g}_{n,\lambda}(z_0)\pm1.96\sigma
\sqrt{(4/3)\cdot I_2/\bigl(n\pi h^\dag\bigr)}.
\end{equation}
Nychka \cite{N88} further shortened (\ref{wahba:bci}) by proposing
%
%e4.4 #&#
\begin{equation}
\label{nychka:bci} \widehat{g}_{n,\lambda}(z_0)\pm1.96\sqrt{\operatorname{Var}
\bigl(b(z_0)\bigr)+\operatorname{Var}\bigl(v(z_0)\bigr)},
\end{equation}
where $b(z_0)=E\{\widehat{g}_{n,\lambda}(z_0)\}-g_0(z_0)$ and
$v(z_0)=\widehat{g}_{n,\lambda}(z_0)-E\{\widehat{g}_{n,\lambda
}(z_0)\}$, and showed that
%
%e4.5 #&#
%e4.6 #&#
\begin{eqnarray}\label{aneq:nychka}
\sigma^2 a\bigl(h^\dag\bigr)/\bigl(\operatorname{Var}
\bigl(b(z_0)\bigr)+\operatorname{Var}\bigl(v(z_0)\bigr)\bigr)
\rightarrow32/27
\\
\eqntext{\mbox{as } n\rightarrow\infty\mbox{ and }\operatorname{Var}\bigl(v(z_0)
\bigr)=8\operatorname{Var}\bigl(b(z_0)\bigr),}
\end{eqnarray}
see his equation (2.3) and the Appendix. Hence, we have
%
%e4.7 #&#
\begin{eqnarray}\label{varres}
\operatorname{Var}\bigl(v(z_0)\bigr)&\sim&\sigma^2\cdot \bigl(I_2/\bigl(n
\pi h^\dag\bigr)\bigr)\quad\mbox{and}
\nonumber\\[-8pt]\\[-8pt]
\operatorname{Var}\bigl(b(z_0)\bigr)&\sim&\bigl(\sigma^2/8\bigr)\cdot
\bigl(I_2/\bigl(n\pi h^\dag\bigr)\bigr).\nonumber
\end{eqnarray}
Therefore, Nychka's Bayesian CI (\ref{nychka:bci}) is asymptotically
equivalent to
%
%e4.8 #&#
\begin{equation}
\label{asympnychka:bci} \widehat{g}_{n,\lambda}(z_0)\pm1.96\sigma
\sqrt{(9/8)\cdot I_2/\bigl(n\pi h^\dag\bigr)}.
\end{equation}

In view of (\ref{asympwahba:bci}) and (\ref{asympnychka:bci}), we
find that the Bayesian CIs of Wahba and Nychka are asymptotically
15.4\% and 6.1\%, respectively, wider than
(\ref{confintformula:refine}). Meanwhile, by~(\ref{varres}) we find
that (\ref{confintformula:refine}) turns out to be a corrected version
of Nychka's CI (\ref{nychka:bci}) by removing the random bias term
$b(z_0)$. The inclusion of $b(z_0)$ in (\ref{nychka:bci}) might be
problematic in that (i) it makes the pointwise limit distribution
nonnormal and thus leads to biased coverage probability; and (ii) it
introduces additional variance, which unnecessarily increases the
length of the interval. By removing $b(z_0)$, we can achieve both
pointwise consistency and a shorter length. Similar considerations
apply when $m>2$. Furthermore, the simulation results in Example
\ref{exampss} demonstrate the superior performance of our CI in both
periodic and nonperiodic splines.

%s4.2 #&#
\subsection{Local likelihood ratio test}\label{seclrt}

In this section, we propose a likelihood ratio method for testing the
value of $g_0(z_0)$ at any $z_0\in\mathbb I$. First, we show that the
null limiting distribution is a scaled noncentral Chi-square with one
degree of freedom. Second, by removing the estimation bias, we obtain a
more useful central Chi-square limit distribution. We also note that as
the smoothness order $m$ approaches infinity, the scaling constant
eventually converges to one. Therefore, we have unveiled an interesting
Wilks phenomenon arising from the proposed nonparametric local testing.
A relevant study was conducted by Banerjee \cite{B07}, who considered a
likelihood ratio test for \textit{monotone} functions, but his
estimation method and null limiting distribution are fundamentally
different from ours.

For some prespecified point $(z_0,w_0)$, we consider the following hypothesis:
%
%e4.9 #&#
\begin{equation}
\label{H0} H_{0}\dvtx g(z_0)=w_0\quad\mbox{versus}\quad H_1\dvtx g(z_0)\neq w_0.
\end{equation}
The ``constrained'' penalized log-likelihood is defined as
$L_{n,\lambda}(g)=\break n^{-1}\sum_{i=1}^n\ell(Y_i; w_0+g(Z_i))-(\lambda/2)
J(g,g)$, where $g\in\mathcal{H}_0=\{g\in\mathcal{H}\dvtx g(z_0)=0\}$. We
consider the likelihood ratio test (LRT) statistic defined as
%
%e4.10 #&#
\begin{equation}\label{LRT:stat}
\mathrm{LRT}_{n,\lambda}=\ell_{n,\lambda}\bigl(w_0+
\widehat{g}{}_{n,\lambda}^{\,0}\bigr)-\ell_{n,\lambda}(
\widehat{g}_{n,\lambda}),
\end{equation}
where $\widehat{g}{}_{n,\lambda}^{\,0}$ is the MLE of $g$ under the local
restriction, that is,
\[
\widehat{g}{}_{n,\lambda}^{\,0}= \arg\max_{g\in
\mathcal{H}_0} L_{n,\lambda}(g).
\]

Endowed with the norm $\|\cdot\|$, $\mathcal{H}_0$ is a closed subset in
$\mathcal{H}$, and thus a Hilbert space. Proposition~\ref{prop:K:W} below
says that $\mathcal{H}_0$ also inherits the reproducing kernel and penalty
operator from $\mathcal{H}$. The proof is trivial and thus omitted.

%
%pr4.2 #&#
\begin{Proposition}\label{prop:K:W}
\textup{(a)} Recall that $K(z_1,z_2)$ is the reproducing kernel for
$\mathcal{H}$ under $\langle\cdot,\cdot\rangle$. The bivariate function
$ K^*(z_1,z_2)=K(z_1,z_2)-(K(z_1,\break z_0)K(z_0,z_2))/ K(z_0,z_0)$ is a
reproducing kernel for $(\mathcal{H}_0,\langle\cdot,\cdot\rangle)$.
That is, for any $z'\in\mathbb{I}$ and $g\in\mathcal{H}_0$, we have $
K^*_{z'}\equiv K^\ast(z',\cdot)\in\mathcal{H}_0$ and $\langle
K_{z'}^*,g\rangle=g(z')$. \textup{(b)} The operator $W_\lambda^*$
defined by $ W_\lambda^*g=W_\lambda g-[(W_\lambda
g)(z_0)/K(z_0,z_0)]\cdot K_{z_0}$ is bounded linear from
$\mathcal{H}_0$ to $\mathcal{H}_0$ and satisfies $\langle W_\lambda^*
g,\widetilde{g}\rangle=\lambda J(g,\widetilde{g})$, for any
$g,\widetilde{g}\in\mathcal{H}_0$.
\end{Proposition}

On the basis of Proposition~\ref{prop:K:W}, we derive the
\textit{restricted} FBR for $\widehat g{}_{n,\lambda}^{\,0}$, which will be
used to obtain the null limiting distribution. By straightforward
calculation we can find the Fr\'{e}chet derivatives of $L_{n,\lambda}$
(under $\mathcal{H}_0$). Let $\Delta g, \Delta g_{j}\in\mathcal{H}_0$
for $j=1,2,3$. The first-order Fr\'{e}chet derivative of
$L_{n,\lambda}$ is
\begin{eqnarray*}
DL_{n,\lambda}(g)\Delta g&=&n^{-1}\sum
_{i=1}^n \dot{\ell}_a
\bigl(Y_i;w_0+g(Z_i)\bigr)\bigl\langle
K_{Z_i}^*,\Delta g\bigr\rangle-\bigl\langle W_\lambda^*g,\Delta
g\bigr\rangle
\\
&\equiv&\bigl\langle S^0_{n}(g), \Delta g\bigr\rangle-
\bigl\langle W_\lambda^*g,\Delta g\bigr\rangle
\\
&\equiv&\bigl\langle
S^0_{n,\lambda}(g), \Delta g\bigr\rangle.
\end{eqnarray*}
Clearly, we have $S^0_{n,\lambda}(\widehat{g}{}_{n,\lambda}^{\,0})=0$.
Define $S^0(g)\Delta g=E\{\langle S_n^0(g),\Delta g\rangle\}$ and\break
$S_\lambda^0(g)\Delta g=S^0(g)\Delta g-\langle W_\lambda^*g,\Delta
g\rangle$. The second-order derivatives are $DS_{n,\lambda}^{\,0}(g)\*\Delta
g_1\Delta g_2=D^2L_{n,\lambda}(g)\Delta g_1\Delta g_2$ and
$DS_\lambda^0(g)\Delta g_1\Delta g_2=DS^0(g)\Delta g_1\Delta
g_2-\break \langle W_\lambda^*\Delta g_1, g_2\rangle$, where
\[
DS^0(g)\Delta g_1\Delta g_2=E\bigl\{\ddot{
\ell}_a\bigl(Y;w_0+g(Z)\bigr)\bigl\langle
K^*_{Z},\Delta g_1\bigr\rangle\bigl\langle
K^*_{Z},\Delta g_2\bigr\rangle\bigr\}.
\]
The third-order Fr\'{e}chet derivative of $L_{n,\lambda}$ is
\begin{eqnarray*}
&& D^3 L_{n,\lambda}(g)\Delta g_1\Delta g_2
\Delta g_3
\\
&&\qquad =n^{-1}\sum_{i=1}^n \ell'''_a
\bigl(Y_i;w_0+g(Z_i)\bigr)\bigl\langle
K^*_{Z_i}, \Delta g_1\bigr\rangle\bigl\langle
K^*_{Z_i},\Delta g_2\bigr\rangle\bigl\langle
K^*_{Z_i},\Delta g_3\bigr\rangle.
\end{eqnarray*}

Similarly to Theorem~\ref{main:thm1}, we need an additional assumption
on the convergence rate of $\widehat{g}{}_{n,\lambda}^{\,0}$:

%asA.4 #&#
\begin{Assumption}\label{A.6}
Under $H_0$,
$\|\widehat{g}{}_{n,\lambda}^{\,0}-g_0^0\|=O_P((nh)^{-1/2}+h^m)$, where
$g_0^0(\cdot)=(g_0(\cdot)-w_0)\in\mathcal{H}_0$.
\end{Assumption}

Assumption~\ref{A.6} is easy to verify by assuming (\ref{A1:a:eq1}),
(\ref{A1:a:eq2}) and $\|\widehat{g}{}_{n,\lambda}^{\,0}-g_0^0\|_{\mathcal
{H}}=o_P(1)$. The proof is similar to that of Proposition
\ref{rates:convergence} by replacing $\mathcal{H}$ with the subspace
$\mathcal{H}_0$.

%
%th4.3 #&#
\begin{Theorem}[(Restricted FBR)]\label{LRT:main:thm1}
Suppose that Assumptions
\ref{A.1},~\ref{A.3}, \ref{A.6} and $H_0$ are satisfied. If $h=o(1)$
and $nh^2\rightarrow\infty$, then
$\|\widehat{g}{}_{n,\lambda}^{\,0}-g_0^0-S_{n,\lambda}^{\,0}(g_0^0)\|
=O_P(a_n\log{n})$.
\end{Theorem}

Our main result on the local LRT is presented below. Define
$r_n=\break(nh)^{-1/2}+h^m$.

%
%th4.4 #&#
\begin{Theorem}[(Local likelihood ratio test)]\label{LRT:lim:dist}
Suppose that Assumptions \mbox{\ref{A.1}--\ref{A.6}} are satisfied. Also
assume $h=o(1)$, $nh^2\rightarrow\infty$,
$a_n=\break o(\min\{r_n, n^{-1}r_n^{-1}(\log{n})^{-1},
n^{-1/2}(\log{n})^{-1}\})$ and $r_n^2h^{-1/2}=o(a_n)$. Furthermore,
for any $z_0\in[0,1]$, $n^{1/2}(W_\lambda
g_0)(z_0)/\sqrt{K(z_0,z_0)}\rightarrow-c_{z_0}$,
%
%e4.11 #&#
\begin{eqnarray}\label{LRT:acondition}
&\displaystyle\lim_{h\rightarrow0}h V(K_{z_0},K_{z_0})\rightarrow
\sigma_{z_0}^2>0\quad\mbox{and}&
\nonumber\\[-8pt]\\[-8pt]
&\displaystyle\lim_{\lambda\rightarrow0}E \bigl\{ I(Z)\bigl|K_{z_0}(Z)\bigr|^2\bigr\}/K(z_0,z_0)
\equiv c_0\in(0,1].& \nonumber
\end{eqnarray}
Under $H_0$, we show: \textup{(i)} $\|\widehat{g}_{n,\lambda}-\widehat
{g}{}_{n,\lambda}^{\,0}-w_0\|=O_P(n^{-1/2})$; \textup{(ii)} $ -2n\cdot
\mathrm{LRT}_{n,\lambda
}=n\|\widehat{g}_{n,\lambda}-\widehat{g}{}_{n,\lambda}^{\,0}-w_0\|
^2+o_P(1)$; and
%
%e4.12 #&#
\begin{equation}\label{main:lim}
\mbox{\textup{(iii)}}\quad {-}2n\cdot
\mathrm{LRT}_{n,\lambda}\stackrel{d} {
\rightarrow}c_0\chi_1^2\bigl(c_{z_0}^2/c_0\bigr)
\end{equation}
with noncentrality parameter $c_{z_0}^2/c_0$.
\end{Theorem}

Note that the parametric convergence rate stated in (i) of
Theorem~\ref{LRT:lim:dist} is reasonable since the restriction is
local. By Proposition~\ref{Prop:RPK}, it can be explicitly
shown that%
%e4.13 #&#
%e4.14 #&#
\begin{eqnarray}\label{finding:c0}
\hspace*{-10pt} c_0=\lim_{\lambda\rightarrow0}
\frac{Q_2(\lambda,z_0)} { Q_1(\lambda,z_0)},
\nonumber\\[-8pt]\\[-8pt]
\eqntext{\mbox{where }
\displaystyle Q_l(\lambda,z)\equiv\sum _{\nu\in\mathbb{N}}
\frac{|h_\nu(z)|^2}{(1+\lambda\gamma_\nu)^l}\qquad\mbox{for }l=1,2.}
\end{eqnarray}
The reproducing kernel $K$, if it exists, is uniquely determined by the
corresponding RKHS; see \cite{D63}. Therefore, $c_0$ defined in
(\ref{LRT:acondition}) depends only on the parameter space. Hence,
different choices of $(\gamma_{\nu}, h_\nu)$ in (\ref{finding:c0}) will
give exactly the same value of $c_0$, although certain choices can
facilitate the calculations. For example, when $\mathcal
H=H_0^m(\mathbb I)$, we can explicitly calculate the value of $c_0$ as
0.75 (0.83) when $m=2$ (3) by choosing the trigonometric polynomial
basis (\ref{H0:basis}). Interestingly, when $\mathcal{H}=H^2(\mathbb
I)$, we can obtain the same value of $c_0$ even without specifying its
(rather different) eigensystem; see Remark~\ref{remark:equivker:appl}
for more details. In contrast, the value of $c_{z_0}$ in
(\ref{main:lim}) depends on the asymptotic bias specified in
(\ref{weta}), whose estimation is notoriously difficult. Fortunately,
under various undersmoothing conditions, we can show $c_{z_0}=0$ and
thus establish a central Chi-square limit distribution. For example, we
can assume higher order smoothness on the true function, as in
Corollary~\ref{LRT:cor2} below.

%
%co4.5 #&#
\begin{Corollary}\label{LRT:cor2}
Suppose that Assumptions~\ref{A.1}--\ref{A.6} are satisfied and $H_0$
holds. Let $m>1+\sqrt{3}/2\approx1.866$. Also assume that the Fourier
coefficients $\{V(g_0,h_\nu)\}_{\nu\in\mathbb{N}}$ of $g_0$ satisfy
$\sum_\nu|V(g_0,h_\nu)|^2 \gamma_\nu^{d}$ for some $d>1+1/(2m)$, which
holds if $g_0\in H^{md}(\mathbb I)$. Furthermore, if
(\ref{LRT:acondition}) is satisfied for any $z_0\in [0,1]$, then
(\ref{main:lim}) holds with the limiting distribution $c_0\chi_1^2$
under $\lambda=\lambda^\ast$.
\end{Corollary}

Corollary~\ref{LRT:cor2} demonstrates a nonparametric type of the Wilks
phenomenon, which approaches the parametric type as
$m\rightarrow\infty$ since $\lim_{m\rightarrow\infty}c_0=1$. This
result provides a theoretical insight for nonparametric local
hypothesis testing; see its \textit{global} counterpart in
Section~\ref{secglobal:lrt:test}.

%s5 #&#
\section{Global asymptotic inference}\label{secjse2}
Depicting the global behavior of a smooth function is crucial in
practice. In Sections~\ref{secsim:confband} and
\ref{secglobal:lrt:test}, we develop the \textit{global} counterparts
of Section~\ref{secjse} by constructing simultaneous confidence bands
and testing a set of global hypotheses.

%s5.1 #&#
\subsection{Simultaneous confidence band}\label{secsim:confband}

In this section, we construct the SCBs for $g$ using the approach of
\cite{BR73}. We demonstrate the theoretical validity of the proposed
SCB based on the FBR and strong approximation techniques. The approach
of \cite{BR73} was originally developed in the context of density
estimation, and it was then extended to M-estimation by \cite{Har89}
and local polynomial estimation by \cite{FZ00,CK03,ZP10}. The
volume-of-tube method \cite{SL94} is another approach, but it
requires the error distribution to be symmetric; see
\cite{ZSW98,KKC10}. Sun et al. \cite{SLM00} relaxed the restrictive
error assumption in generalized linear models, but they had to
translate the nonparametric estimation into parametric estimation. Our
SCBs work for a general class of nonparametric models including normal
regression and logistic regression. Additionally, the minimum width of
the proposed SCB is shown to achieve the lower bound established by
\cite{GW08}; see Remark~\ref{scb:minimax}. An interesting by-product is
that, under the equivalent kernel conditions given in this section, the
local asymptotic inference procedures developed from cubic splines and
periodic splines are essentially the same despite the intrinsic
difference in their eigensystems; see Remark~\ref{remark:equivker:appl}
for technical details.

The key conditions assumed in this section are the equivalent kernel
conditions (\ref{equiv:kernel})--(\ref{covfun}). Specifically, we
assume that there exists a real-valued function $\omega(\cdot)$ defined
on $\mathbb{R}$ satisfying, for any fixed $0<\varphi<1$, $h^\varphi
\leq z\leq1- h^\varphi$ and $t\in\mathbb{I}$,
%
%e5.1 #&#
\begin{eqnarray}\label{equiv:kernel}
&& \biggl\llvert\frac{d^j}{dt^j}\bigl(h^{-1} \omega\bigl((z-t)/ h\bigr)-K(z,t) \bigr)\biggr\rrvert
\nonumber\\[-8pt]\\[-8pt]
&&\qquad \leq
C_K h^{-(j+1)}\exp\bigl(-C_2 h^{-1+\varphi}\bigr)\qquad\mbox{for } j=0,1,\nonumber
\end{eqnarray}
where $C_2,C_K$ are positive constants. Condition (\ref{equiv:kernel})
implies that $\omega$ is an equivalent kernel of the reproducing
kernel function $K$ with a certain degree of approximation accuracy. We
also require two regularity conditions on $\omega$:
%
%e5.2 #&#
%e5.3 #&#
\begin{eqnarray}\label{convfun2}
\bigl|\omega(u)\bigr|\leq C_\omega\exp\bigl(-|u|/C_3\bigr),\qquad
\bigl|\omega'(u)\bigr|\leq C_\omega\exp\bigl(-|u|/C_3\bigr)
\nonumber\\[-8pt]\\[-8pt]
\eqntext{\mbox{for any } u\in\mathbb{R},}
\end{eqnarray}
and there exists a constant $0<\rho\leq2$ s.t.
%
%e5.4 #&#
\begin{equation}\label{covfun}
\int_{-\infty}^\infty\omega(t)\omega(t+z)\,dt=
\sigma_\omega^2-C_\rho|z|^\rho+o
\bigl(|z|^\rho\bigr)\qquad\mbox{as }|z|\rightarrow\infty,
\end{equation}
where $C_3, C_\omega, C_\rho$ are positive constants and $\sigma
_\omega^2=\int_\mathbb{R}\omega(t)^2\,dt$. An example of $\omega$
satisfying (\ref{equiv:kernel})--(\ref{covfun}) will be given in
Proposition~\ref{prop:equiv:ker}. The following exponential envelope
condition is also needed:
%
%e5.5 #&#
\begin{equation}\label{envelop:exp}
\sup_{z,t\in\mathbb{I}} \biggl|\frac{\partial}{\partial
z}K(z,t) \biggr|=O\bigl(h^{-2}\bigr).
\end{equation}

%
%th5.1 #&#
\begin{Theorem}[(Simultaneous confidence band)]\label{confband}
Suppose Assumptions \mbox{\ref{A.1}--\ref{A.5}} are satisfied, and $Z$ is
uniform on $\mathbb{I}$. Let $m>(3+\sqrt{5})/4\approx 1.3091$ and
$h=n^{-\delta}$ for any $\delta\in(0,2m/(8m-1))$. Furthermore,
$E\{\exp(|\epsilon|/C_0)\mid Z\}\leq C_1$, a.s., and
(\ref{equiv:kernel})--(\ref{envelop:exp}) hold. The conditional density
of $\epsilon$ given $Z=z$, denoted $\pi (\epsilon\mid z)$, satisfies the
following: for some positive constants $\rho_1$ and~$\rho_2$,
%
%e5.6 #&#
\begin{equation}\label{conddensity}
\biggl\llvert\frac{d}{dz}\log{\pi(\epsilon
\mid z)}\biggr
\rrvert\leq\rho_1\bigl(1+|\epsilon|^{\rho_2}\bigr)\qquad\mbox{for
any } \epsilon\in\mathbb{R}\mbox{ and } z\in\mathbb{I}.
\end{equation}
Then, for any $0<\varphi<1$ and $u\in\mathbb{R}$,
%
%e5.7 #&#
\begin{eqnarray}\label{ssband:eq}
\qquad&&P \Bigl((2\delta\log{n})^{1/2} \Bigl\{\sup
_{h^\varphi\leq z\leq
1-h^\varphi}(nh)^{1/2}\sigma_\omega^{-1}I(z)^{-1/2}\nonumber
\\
&&\hspace*{117pt}{}\times \bigl|
\widehat{g}_{n,\lambda}(z)-g_0(z)+(W_\lambda
g_0) (z)\bigr|-d_n \Bigr\}\leq u \Bigr)
\\
&&\qquad \longrightarrow\exp\bigl(-2\exp(-u)\bigr),\nonumber
\end{eqnarray}
where $d_n$ relies only on $h$, $\rho$, $\varphi$ and $C_\rho$.
\end{Theorem}

The FBR developed in Section~\ref{secfbr} and the strong approximation
techniques developed by \cite{BR73} are crucial to the proof of
Theorem~\ref{confband}. The uniform distribution on $Z$ is assumed only
for simplicity, and this condition can be relaxed by requiring that the
density is bounded away from zero and infinity. Condition
(\ref{conddensity}) holds in various situations such as the conditional
normal model $\epsilon\mid Z=z\sim N(0,\sigma^2(z))$, where
$\sigma^2(z)$ satisfies $\inf_z \sigma^2(z)>0$, and $\sigma(z)$ and
$\sigma'(z)$ both have finite upper bounds. The existence of the bias
term $W_\lambda g_0(z)$ in (\ref{ssband:eq}) may result in poor
finite-sample performance. We address this issue by assuming
undersmoothing, which is advocated by \cite{NP98,Hall91,Hall92}; they
showed that undersmoothing is more efficient than explicit bias
correction when the goal is to minimize the coverage error.
Specifically, the bias term will asymptotically vanish if we assume
that
%
%e5.8 #&#
\begin{equation}\label{confband:removebias} \lim_{n\rightarrow\infty} \Bigl\{\sup
_{h^\varphi\leq z\leq
1-h^\varphi}\sqrt{nh\log n}\bigl|W_\lambda g_0(z)\bigr|
\Bigr\}=0.
\end{equation}
Condition~(\ref{confband:removebias}) is slightly stronger than the
undersmoothing condition\break $\sqrt{nh}(W_\lambda g_0)(z_0)=o(1)$ assumed
in Proposition~\ref{Prop:conf}. By the generalized Fourier expansion of
$W_\lambda g_0$ and the uniform boundedness of $h_\nu$ (see
Assumption~\ref{A.3}), we can show that (\ref{confband:removebias})
holds if we properly increase the amount of smoothness on $g_0$ or
choose a suboptimal $\lambda$, as in Corollaries
\ref{l2:regression:CLT}~and~\ref{LRT:cor2}.

Proposition~\ref{prop:equiv:ker} below demonstrates the validity of
Conditions (\ref{equiv:kernel})--(\ref{covfun}) in $L_2$ regression.
The proof relies on an explicit construction of the equivalent kernel
function obtained by \cite{MG93}. We consider only $m=2$ for simplicity.

%
%pr5.2 #&#
\begin{Proposition}[($L_2$ regression)]\label{prop:equiv:ker}
Let $\ell(y;a)=-(y-a)^2/(2\sigma^2)$, $Z\sim
\operatorname{Unif}[0,1]$ and $\mathcal H=H^2(\mathbb I)$, that is,
$m=2$. Then, (\ref{equiv:kernel})--(\ref{covfun}) hold with
$\omega(t)=\sigma ^{2-1/m}\omega_0(\sigma^{-1/m}t)$ for
$t\in\mathbb{R}$, where\vspace*{-1pt}
$\omega_0(t)=\frac{1}{2\sqrt{2}}\exp(-|t|/\sqrt{2}) (\cos(t/\sqrt
{2})+\sin(|t|/\sqrt{2}) )$. In particular, (\ref{covfun}) holds for
arbitrary $\rho\in(0,2]$ and $C_\rho=0$.
\end{Proposition}

%
%re5.1 #&#
\begin{Remark}\label{remark:scb}
In the setting of Proposition~\ref{prop:equiv:ker}, we are able to
explicitly find the constants $\sigma_\omega^2$ and $d_n$ in
Theorem~\ref{confband}. Specifically, by direct calculation it can~be
found that $\sigma_\omega^2=0.265165\sigma^{7/2}$ since $\sigma_{\omega
_0}^2=\int_{-\infty}^\infty|\omega_0(t)|^2\,dt=0.265165$ when $m=2$.
Choose $C_\rho=0$ for arbitrary $\rho\in(0,2]$. By the formula of
$B(t)$ given in Theorem~A.1 of \cite{BR73}, we know that
%
%e5.9 #&#
\begin{equation}\label{l2dn}
\quad d_n=\bigl(2\log\bigl(h^{-1}-2h^{\varphi-1}
\bigr)\bigr)^{1/2}+\frac{(1/\rho-1/2)
\log\log(h^{-1}-2h^{\varphi-1})}{(2\log(h^{-1}-2h^{\varphi-1}))^{1/2}}.
\end{equation}
When $\rho=2$, the above $d_n$ is simplified as $(2\log
(h^{-1}-2h^{\varphi-1}))^{1/2}$. In general, we observe that
$d_n\sim(-2\log{h})^{1/2}\asymp\sqrt{\log n}$ for sufficiently
large $n$ since $h=n^{-\delta}$. Given that the estimation bias is
removed, for example, under (\ref{confband:removebias}), we obtain the
following $100\times(1-\alpha)\%$ SCB:
%
%e5.10 #&#
\begin{eqnarray}\label{ss:confband}
&& \bigl\{ \bigl[\widehat{g}_{n,\lambda}(z)\pm
0.5149418(nh)^{-1/2}\widehat{\sigma}{}^{3/4} \bigl(c_\alpha^*/
\sqrt{-2\log{h}}+d_n \bigr) \bigr]\dvtx
\nonumber\\[-8pt]\\[-8pt]
&& \hspace*{187pt} h^\varphi\leq z \leq1-h^\varphi\bigr\},\nonumber
\end{eqnarray}
where $d_n=(-2\log{h})^{1/2}$, $c_\alpha^*=-\log(-\log(1-\alpha )/2)$
and $\widehat{\sigma}$ is a consistent estimate of $\sigma$. Therefore,
to obtain uniform coverage, we have to increase the bandwidth up to an
order of $\sqrt{\log{n}}$ over the length of the pointwise CI given in
(\ref{confintformula:refine}). Note that we have excluded the boundary
points in (\ref{ss:confband}).
\end{Remark}

%
%re5.2 #&#
\begin{Remark}\label{remark:equivker:appl}
An interesting by-product we discover in the setting of Proposition
\ref{prop:equiv:ker} is that the pointwise asymptotic CIs for
$g_0(z_0)$ based on cubic splines and periodic splines share the same
length at any $z_0\in(0,1)$. This result is surprising since the two
splines have intrinsically different structures. Under
(\ref{equiv:kernel}), it can be shown that
\begin{eqnarray*}
\sigma_{z_0}^2&\sim&\sigma^{-2} h\int_0^1 \bigl|K(z_0,z)\bigr|^2\,dz
\\
&\sim&\sigma^{-2}h^{-1} \int_0^1\biggl\llvert\omega\biggl(\frac
{z-z_0}{h} \biggr)\biggr\rrvert^2\,dz
\\
&=&\sigma^{-2}\int_{-z_0/h}^{(1-z_0)/h}\bigl|
\omega(s)\bigr|^2\,ds\sim\sigma^{-2}\int_\mathbb{R}\bigl|
\omega(s)\bigr|^2\,ds=\sigma^{3/2}\sigma_{\omega_0}^2,
\end{eqnarray*}
given the choice of $\omega$ in Proposition~\ref{prop:equiv:ker}.
Thus, Corollary~\ref{main:thm3} implies the following 95\% CI:
%
%e5.11 #&#
\begin{equation}
\label{cubic:confint}
\qquad \widehat{g}_{n,\lambda}(z_0)\pm1.96
(nh)^{-1/2}\sigma^{3/4} \sigma_{\omega_0} =
\widehat{g}_{n,\lambda}(z_0)\pm1.96 \bigl(nh^\dag
\bigr)^{-1/2}\sigma\sigma_{\omega_0}.
\end{equation}
Since $\sigma_{\omega_0}^2=I_2/\pi$, the lengths of the CIs
(\ref{confintformula:refine}) (periodic spline) and
(\ref{cubic:confint}) (cubic spline) coincide with each other. The
above calculation of $\sigma_{z_0}^2$ relies on $L_2$ regression. For
general models such as logistic regression, one can instead use
a~weighted version of (\ref{est:f}) with the weights $B(Z_i)^{-1}$ to
obtain the exact formula. Another application of Proposition
\ref{prop:equiv:ker} is to find the value of $c_0$ in
Theorem~\ref{LRT:main:thm1} for the construction of the local LRT test
when $\mathcal H=H^2(\mathbb I)$. According to the definition of $c_0$,
that is, (\ref{LRT:acondition}), we have
$c_0\sim\sigma_{z_0}^2/(hK(z_0,z_0))$. \mbox{Under} (\ref{equiv:kernel}), we
have $K(z_0,z_0)\sim h^{-1}\omega(0)=h^{-1}\sigma^{3/2}\omega
_0(0)=0.3535534h^{-1}\sigma^{3/2}$. Since
$\sigma_{z_0}^2\sim\sigma^{3/2}\sigma_{\omega_0}^2$ and
$\sigma_{\omega_0}^2=I_2/\pi$, we have $c_0= 0.75$. This value
coincides with that found in periodic splines, that is, $\mathcal
H=H^2_0(\mathbb I)$. These intriguing phenomena have never been
observed in the literature and may be useful for simplifying the
construction of CIs and local LRT.
\end{Remark}

%
%re5.3 #&#
\begin{Remark}\label{scb:minimax}
Genovese and Wasserman \cite{GW08} showed that when $g_0$ belongs to an
$m$th-order Sobolev ball, the lower bound for the average width of an
SCB is proportional to $b_n n^{-m/(2m+1)}$, where $b_n$ depends only on
$\log n$. We next show that the (minimum) bandwidth of the proposed SCB
can achieve this lower bound with $b_n=(\log{n})^{(m+1)/(2m+1)}$. Based
on Theorem~\ref{confband}, the width of the SCB is of order $d_n
(nh)^{-1/2}$, where $d_n\asymp\sqrt{\log{n}}$; see Remark
\ref{remark:scb}. Meanwhile, Condition (\ref{confband:removebias}) is
crucial for our band to maintain the desired coverage probability.
Suppose that the Fourier coefficients of $g_0$ satisfy $\sum_\nu
|V(g_0,h_\nu)|\gamma_\nu^{1/2}<\infty$. It can be verified that
(\ref{confband:removebias}) holds when $nh^{2m+1}\log{n}=O(1)$, which
sets an upper bound for $h$, that is,\break $O(n\log{n})^{-1/(2m+1)}$. When
$h$ is chosen as the above upper bound and $d_n\asymp\sqrt{\log n}$,
our SCB achieves the minimum order of bandwidth\break $n^{-m/(2m+1)}(\log
{n})^{(m+1)/(2m+1)}$, which turns out to be optimal according to~\cite{GW08}.
\end{Remark}

In practice, the construction of our SCB requires a delicate choice of
$(h, \varphi)$. Otherwise, over-coverage or undercoverage of the true
function may occur near the boundary points. There is no practical or
theoretical guideline on how to find the optimal $(h,\varphi)$,
although, as noted by \cite{BR73}, one can choose a proper $h$ to make
the band as thin as possible. Hence, in the next section, we propose a
more straightforward
likelihood-ratio-based approach for testing the global behavior, which
requires only tuning $h$.

%s5.2 #&#
\subsection{Global likelihood ratio test}\label{secglobal:lrt:test}

There is a vast literature dealing with nonparametric hypothesis
testing, among which the GLRT proposed by Fan et al.~\cite{FZZ01}
stands out. Because of the technical complexity, they focused on the
local polynomial fitting; see \cite{FZ04} for a sieve version. Based on
smoothing spline estimation, we propose the PLRT, which is applicable
to both simple and composite hypotheses. The null limiting distribution
is identified to be nearly Chi-square with diverging degrees of
freedom. The degrees of freedom depend only on the functional parameter
space, while the null limiting distribution of the GLRT depends on the
choice of kernel functions; see Table~2 in \cite{FZZ01}. Furthermore,
the PLRT is shown to achieve the minimax rate of testing in the sense
of \cite{I93}. As demonstrated in our simulations, the PLRT performs
better than the GLRT in terms of power, especially in small-sample
situations. Other smoothing-spline-based testing such as LMP, GCV and
GML (see \cite{CKWY88,W90,J96,C94,RG00,LW02}) use ad-hoc discrepancy
measures leading to complicated null distributions involving nuisance
parameters; see a thorough review in \cite{LW02}.

Consider the following ``global'' hypothesis:
%
%e5.12 #&#
\begin{equation}\label{global:H0}
H_{0}^{\mathrm{global}}\dvtx g=g_0\quad\mbox{versus}\quad
H_1^{\mathrm{global}}\dvtx g\in\mathcal H-\{g_0\},
\end{equation}
where $g_0\in\mathcal H$ can be either known or unknown. The PLRT
statistic is defined to be
%
%e5.13 #&#
\begin{equation}
\label{Global:LRT:Test} \mathrm{PLRT}_{n,\lambda}=\ell_{n,\lambda}(g_0)-
\ell_{n,\lambda}(\widehat{g}_{n,\lambda}).
\end{equation}

Theorem~\ref{Global:LRT:Lim} below derives the null limiting
distribution of $\mathrm{PLRT}_{n,\lambda}$. We remark that the null
limiting distribution remains the same even when the hypothesized value
$g_0$ is unknown (whether its dimension is finite or infinite). This
nice property can be used to test the composite hypothesis; see
Remark~\ref{rem:comp}.

%
%th5.3 #&#
\begin{Theorem}[(Penalized likelihood ratio test)]\label{Global:LRT:Lim}
Let Assumptions~\ref{A.1}--\ref{A.5} be~satisfied. Also assume
$nh^{2m+1}=O(1)$, $nh^2\rightarrow\infty$, $a_n=o(\min\{
r_n,\break n^{-1}r_n^{-1}\* h^{-1/2}(\log{n})^{-1}, n^{-1/2}(\log{n})^{-1}\})$
and $r_n^2h^{-1/2}=o(a_n)$. Furthermore, under $H_0^{\mathrm{global}}$,
$E\{\epsilon^4\mid Z\}\leq C$, a.s., for some constant $C>0$, where
$\epsilon=\dot{\ell}_a(Y;g_0(Z))$ represents the ``model error.'' Under
$H_0^{\mathrm{global}}$, we have
%
%e5.14 #&#
\begin{equation}
\label{lim:dist:globalLRT} (2u_n)^{-1/2}
\bigl(-2nr_K\cdot \mathrm{PLRT}_{n,\lambda}-nr_K
\|W_\lambda g_0\| ^2-u_n \bigr)
\stackrel{d} {\longrightarrow}N(0,1),
\end{equation}
where $u_n=h^{-1}\sigma_K^4/\rho_K^2$, $r_K=\sigma_K^2/\rho_K^2$,
%
%e5.15 #&#
\begin{eqnarray}\label{sigma:rho:formula}
\sigma_K^2&=& hE\bigl\{
\epsilon^2 K(Z,Z)\bigr\}=\sum_\nu \frac
{h}{(1+\lambda\gamma_\nu)},
\nonumber\\[-8pt]\\[-8pt]
\rho_K^2&=&hE\bigl\{ \epsilon_1^2\epsilon_2^2
K(Z_1,Z_2)^2\bigr\}=\sum _\nu\frac{h}{(1+\lambda\gamma_\nu)^2}\nonumber
\end{eqnarray}
and $(\epsilon_i,Z_i)$, $i=1,2$ are i.i.d. copies of $(\epsilon,Z)$.
\end{Theorem}

A direct examination reveals that $h\asymp n^{-d}$ with $\frac
{1}{2m+1}\leq d<\frac{2m}{8m-1}$ satisfies the rate conditions required
by Theorem~\ref{Global:LRT:Lim} when $m>(3+\sqrt{5})/4\approx1.309$. By
the proof of Theorem~\ref{Global:LRT:Lim}, it can be shown that
$n\|W_\lambda g_0\|^2=o(h^{-1})=o(u_n)$. Therefore, $-2nr_K\cdot
\mathrm{PLRT}_{n,\lambda}$ is asymptotically $N(u_n,2u_n)$. As $n$
approaches infinity, $N(u_n,2u_n)$ is nearly $\chi^2_{u_n}$. Hence,
$-2nr_K\cdot \mathrm{PLRT}_{n,\lambda}$ is approximately distributed as
$\chi_{u_n}^2$, denoted
%
%e5.16 #&#
\begin{equation}
\label{eq:global:LRT} -2nr_K\cdot \mathrm{PLRT}_{n,\lambda}\stackrel{a} {
\sim} \chi^2_{u_n}.
\end{equation}
That is, the Wilks phenomenon holds for the PLRT. The specifications of
(\ref{eq:global:LRT}), that is, $\sigma_K^2$ and $\rho_K^2$, are
determined only by the parameter space and model setup. We also note
that undersmoothing is not required for our global test.

We next discuss the calculation of $(r_K, u_n)$. In the setting of
Proposition~\ref{prop:equiv:ker}, it can be shown by the equivalent
kernel conditions that $\sigma_K^2=h\sigma^{-2} \int_0^1 K(z,z)\,dz\sim
h\sigma^{-2} (h^{-1}\omega(0))= \sigma^{-1/2} \omega
_0(0)=0.3535534\sigma^{-1/2}$ and\break $\rho_K^2\sim\sigma^{-1/2}\sigma
_{\omega_0}^2= 0.265165\sigma^{-1/2}$. So $r_K=1.3333$ and $u_n=0.4714
h^{-1}\sigma^{-1/2}$. If we replace $H^2(\mathbb I)$ by $H_0^2(\mathbb
I)$, direct calculation in case~(I) of Example~\ref{exampss} reveals
that $(r_K, u_n)$ have exactly the same values. When $\mathcal
H=H^m_0(\mathbb I)$, we have $2r_K\rightarrow2$ as $m$ tends to
infinity. This limit is consistent with the scaling constant two in the
parametric likelihood ratio theory. In $L_2$ regression, the possibly
unknown parameter $\sigma$ in $u_n$ can be profiled out without
changing the null limiting distribution. In practice, by the wild
bootstrap we can directly simulate the null limiting distribution by
fixing the nuisance parameters at some reasonable values or estimates
without finding the values of $(r_K, u_n)$. This is a major advantage
of the Wilks type of results.

%
%re5.4 #&#
\begin{Remark}\label{rem:comp}
We discuss composite hypothesis testing via the PLRT. Specifically, we
test whether $g$ belongs to some finite-dimensional class of functions,
which is much larger than the null space $\mathcal N_m$ considered in
the literature. For instance, for any integer $q\geq0$, consider the
null hypothesis
%
%e5.17 #&#
\begin{equation}
\label{compo} H_0^{\mathrm{global}}\dvtx g\in\mathcal{L}_q(
\mathbb{I}),
\end{equation}
where $\mathcal{L}_q(\mathbb{I})\equiv\{g(z)=\sum_{l=0}^q a_l z^l\dvtx
a=(a_0,a_1,\ldots,a_q)^T\in\mathbb{R}^{q+1}\}$ is the class of~the
$q$th-order polynomials. Let
$\widehat{a}_*=\arg\max_{a\in\mathbb{R}^{q+1}}
\{(1/n)\sum_{i=1}^n\ell(Y_i;\break \sum_{l=0}^q a_l Z_i^l)-(\lambda/2) a^T D
a\}$, where
\[
D=\int_0^1 \bigl(0,0, 2, 6z, \ldots,
q(q-1)z^{q-2}\bigr)^T \bigl(0,0, 2, 6z, \ldots,
q(q-1)z^{q-2}\bigr)\,dz
\]
is a $(q+1)\times(q+1)$ matrix. Hence, under $H_0^{\mathrm{global}}$,
the penalized MLE is $\widehat{g}_*(z)=\sum_{l=0}^q \widehat{a}_{*l}
z^l$. Let $g_{0q}$ be an unknown ``true parameter'' in
$\mathcal{L}_q(\mathbb {I})$ corresponding to a vector of polynomial
coefficients $a^0=(a^0_0,a^0_1,\ldots,a^0_q)^T$. To test~(\ref{compo}),
we decompose the PLRT statistic as
$\mathrm{PLRT}_{n,\lambda}^{\mathrm{com}}=L_{n1}-L_{n2}$, where
$L_{n1}=\ell_{n,\lambda}(g_{0q})-\ell_{n,\lambda}(\widehat
{g}_{n,\lambda})$ and
$L_{n2}=\ell_{n,\lambda}(g_{0q})-\ell_{n,\lambda}(\widehat{g}_*)$. When
we formulate
\[
H_0'\dvtx a=a^0\quad\mbox{versus}\quad H_1'\dvtx a\neq a^0,
\]
$L_{n2}$ appears to be the PLRT statistic in the parametric setup. It
can be shown that $L_{n2}=O_P(n^{-1})$ whether $q<m$ (by applying the
parametric theory in \cite{Shao03}) or $q\geq m$ (by slightly
modifying the proof of Theorem~\ref{LRT:lim:dist}). On the other hand,
$L_{n1}$~is exactly the PLRT for testing
\[
H_0'\dvtx g=g_{0q}\quad\mbox{versus}\quad H_1^{\mathrm{global}}\dvtx g\neq g_{0q}.
\]
By Theorem~\ref{Global:LRT:Lim}, $L_{n1}$ follows the limit
distribution specified in (\ref{eq:global:LRT}). In summary, under
(\ref{compo}), $\mathrm{PLRT}_{n,\lambda}^{\mathrm{com}}$ has the same limit
distribution since $L_{n2}=O_P(n^{-1})$ is negligible.
\end{Remark}

To conclude this section, we show that the PLRT achieves the optimal
minimax rate of testing specified in Ingster \cite{I93} based on a
\textit{uniform} version of the FBR. For convenience, we consider only
$\ell(Y;a)=-(Y-a)^2/2$. Extensions to a more general setup can be found
in the supplementary document \cite{SC12} under stronger assumptions,
for example, a more restrictive alternative set.

Write the local alternative as $H_{1n}\dvtx g=g_{n0}$, where
$g_{n0}=g_0+g_n$, \mbox{$g_0\in\mathcal{H}$} and $g_{n}$ belongs to the
alternative value set $\mathcal{G}_a\equiv\{g\in\mathcal{H}\mid
\operatorname{Var}(g(Z)^2)\leq\zeta E^2\{g(Z)^2\}, J(g,g)\leq\zeta\}$
for some constant $\zeta>0$.

%
%th5.4 #&#
\begin{Theorem}\label{power:thm}
Let\vspace*{-1pt} $m>(3+\sqrt {5})/4\approx1.309$ and $h\asymp n^{-d}$ for
$\frac{1}{2m+1}\leq d<\frac{2m}{8m-1}$. Suppose that Assumption
\ref{A.3} is satisfied, and uniformly over $g_{n}\in\mathcal{G}_a$,
$\|\widehat {g}_{n,\lambda}-g_{n0}\|=O_\Pr(r_n)$ holds under
$H_{1n}\dvtx g=g_{n0}$. Then for any $\delta\in(0,1)$, there exist
positive constants $C$ and $N$ such that
%
%e5.18 #&#
\begin{equation}
\label{power:inequality} \inf_{n\geq N} \mathop{\inf
_{g_{n}\in\mathcal{G}_a}}_{\|g_n\|\geq C\eta_n} \Pr\bigl(\mbox{reject }
H_0^{\mathrm{global}}\mid  \mbox{$H_{1n}$ is true} \bigr)\geq1-
\delta,
\end{equation}
where $\eta_n\geq\sqrt{h^{2m}+(nh^{1/2})^{-1}}$. The minimal lower
bound of $\eta_n$, that is, $n^{-2m/(4m+1)}$, is achieved when
$h=h^{\ast\ast}\equiv n^{-2/(4m+1)}$.
\end{Theorem}

The condition ``uniformly over $g_{n}\in\mathcal{G}_a$,
$\|\widehat{g}_{n,\lambda}-g_{n0}\|=O_\Pr(r_n)$ holds under $H_{1n}\dvtx
g=g_{n0}$'' means that for any $\widetilde\delta>0$, there exist
constants $\widetilde C$ and $\widetilde N$, both unrelated to
$g_{n}\in\mathcal G_a$, such that $\inf_{n\geq\widetilde
N}\inf_{g_{n}\in\mathcal {G}_a}P_{g_{n0}}
(\|\widehat{g}_{n,\lambda}-g_{n0}\|\leq\widetilde Cr_n
)\geq1-\widetilde\delta$.

Theorem~\ref{power:thm} proves that, when $h=h^{\ast\ast}$, the PLRT
can detect any local alternatives with separation rates no faster than
$n^{-2m/(4m+1)}$, which turns out to be the minimax rate of testing in
the sense of Ingster \cite{I93}; see Remark~\ref{rem:mini} below.

%
%re5.5 #&#
\begin{Remark}\label{rem:mini}
The minimax rate of testing established in Ingster \cite{I93} is
\mbox{under} the usual $\|\cdot\|_{L_2}$-norm (w.r.t. the Lebesgue measure).
However, the separation rate derived under the $\|\cdot\|$-norm is
still optimal because of the trivial domination of $\|\cdot\|$ over $\|
\cdot\|_{L_2}$ (under the conditions of Theorem~\ref{power:thm}). Next,
we heuristically explain why the minimax rates of testing associated
with $\|\cdot\|$, denoted $b_n'$, and with $\|\cdot\| _{L_2}$, denoted
$b_n$, are the same. By definition, whenever $\|g_n\| \geq b_n'$ or
$\|g_n\|_{L_2}\geq b_n$, $H_0^{\mathrm{global}}$ can be rejected with a
large probability, or equivalently, the local alternatives can be
detected. $b_n'$ and $b_n$ are the minimum rates that satisfy this
property. Ingster \cite{I93} has shown that $b_n\asymp n^{-2m/(4m+1)}$.
Since $\|g_n\|_{L_2}\geq b_n'$ implies $\|g_n\|\geq b_n'$,
$H_0^{\mathrm{global}}$ is rejected. This means $b_n'$ is an upper
bound for detecting the local alternatives in terms of
$\|\cdot\|_{L_2}$ and so $b_n\leq b_n'$. On the other hand, suppose
$h=h^{\ast\ast}\asymp n^{-2/(4m+1)}$ and $\|g_n\|\geq C
n^{-2m/(4m+1)}\asymp b_n$ for some large $C>\zeta^{1/2}$. Since
$\lambda J(g_n,g_n)\leq\zeta\lambda\asymp\zeta n^{-4m/(4m+1)}$, it
follows that $\|g_n\|_{L_2}\geq(C^2-\zeta)^{1/2} n^{-2m/(4m+1)}\asymp
b_n$. This means $b_n$ is a upper bound for detecting the local
alternatives in terms of $\|\cdot\|$ and so $b_n'\leq b_n$. Therefore,
$b_n'$ and $b_n$ are of the same order.
\end{Remark}

%s6 #&#
\section{Examples}\label{secexa}

In this section, we provide three concrete examples together with simulations.

%
%ex6.1 #&#
\begin{Example}[($L_2$ regression)]\label{exampss}
We consider the regression model with an additive error
%
%e6.1 #&#
\begin{equation}
\label{npmodel} Y=g_0(Z)+\epsilon,
\end{equation}
where $\epsilon\sim N(0,\sigma^2)$ with an unknown variance $\sigma
^2$. Hence, $I(Z)=\sigma^{-2}$ and $V(g, \widetilde
g)=\sigma^{-2}E\{g(Z)\widetilde g(Z)\}$. For simplicity, $Z$ was
generated uniformly over $\mathbb{I}$. The function \textit{ssr}(\,) in the R
package \textit{assist} was used to select the smoothing parameter
$\lambda$ based on CV or GCV; see \cite{KW02}. We first consider
$\mathcal H=H_0^m(\mathbb{I})$ in case~(I) and then $\mathcal
H=H^m(\mathbb{I})$ in case~(II).
\begin{longlist}[Case (I)]
\item[Case (I). $\mathcal H=H_0^m(\mathbb{I})$:] In this
    case, we choose the basis functions as
%
%e6.2 #&#
\begin{equation}\label{H0:basis}
h_\mu(z)=\cases{
\sigma,&\quad $\mu=0$,
\vspace*{2pt}\cr
\sqrt{2}\sigma\cos(2\pi kz), &\quad $\mu=2k, k=1,2,\ldots,$
\vspace*{2pt}\cr
\sqrt{2}\sigma\sin(2\pi kz), &\quad $\mu=2k-1, k=1,2,\ldots,$}
\end{equation}
with the corresponding eigenvalues $\gamma_{2k-1}=\gamma_{2k}=\sigma
^2(2\pi k)^{2m}$ for $k\geq1$ and \mbox{$\gamma_0=0$}. Assumption~\ref{A.3}
trivially holds for this choice of $(h_\mu, \gamma_\mu)$. The lemma
below is useful for identifying the critical quantities for inference.
\end{longlist}

%
%le6.1 #&#
\begin{Lemma}\label{lemeg1}
Let $I_l=\int_0^\infty(1+x^{2m})^{-l}\,dx$ for $l=1,2$ and $h^\dag
=h\sigma^{1/m}$. Then
%
%e6.3 #&#
\begin{equation}\label{pss:aeq}
\sum_{k=1}^{\infty}\frac{1}{(1+(2\pi
h^\dag k)^{2m})^l}\sim\frac{I_l}{2\pi h^\dag}.
\end{equation}
\end{Lemma}

By Proposition~\ref{Prop:conf}, the asymptotic 95\% pointwise CI for
$g(z_0)$ is $\widehat{g}_{n,\lambda}(z_0)\pm1.96\sigma_{z_0}/\sqrt{nh}$
when ignoring the bias. By the definition of $\sigma_{z_0}^2$ and
Lemma~\ref{lemeg1}, we have
\[
\sigma_{z_0}^2\sim hV(K_{z_0},K_{z_0})=
\sigma^2 h \Biggl(1+2 \sum_{k=1}^{\infty}
\bigl(1+\bigl(2\pi h^\dag k\bigr)^{2m}\bigr)^{-2}
\Biggr)\sim\bigl(I_2\sigma^{2-1/m}\bigr)/\pi.
\]
Hence, the CI becomes
%
%e6.4 #&#
\begin{equation}\label{confintformula}
\widehat{g}_{n,\lambda}(z_0)\pm1.96
\widehat\sigma{}^{1-1/(2m)}\sqrt{I_2/(\pi nh)},
\end{equation}
where $\widehat{\sigma}{}^2=\sum_i
(Y_i-\widehat{g}_{n,\lambda}(Z_i))^2/(n-\operatorname{trace}(A(\lambda)))$ is a
consistent estimate of $\sigma^2$ and $A(\lambda)$ denotes the
smoothing matrix; see \cite{W90}. By (\ref{finding:c0}) and
(\ref{H0:basis}), for $l=1,2$,
\begin{eqnarray*}
Q_l(\lambda,z_0)&=&\sigma^2+\sum
_{k\geq1} \biggl\{\frac
{|h_{2k}(z_0)|^2}{(1+\lambda\sigma^2(2\pi
k)^{2m})^l}+\frac{|h_{2k-1}(z_0)|^2}{(1+\lambda\sigma^2(2\pi
k)^{2m})^l} \biggr\}
\\
&=&\sigma^2+2\sigma^2\sum_{k\geq1}
\frac{1}{(1+\lambda\sigma
^2(2\pi k)^{2m})^l}
\\
&=&\sigma^2+2\sigma^2\sum_{k\geq1} \frac
{1}{(1+(2\pi h^\dag k)^{2m})^l}.
\end{eqnarray*}
By Lemma~\ref{lemeg1}, we have $c_0=I_2/I_1$. In particular, $c_0=0.75$
(0.83) when \mbox{$m=2$ (3)}.
%
%f1 #&#
\begin{figure}%[b]

\includegraphics{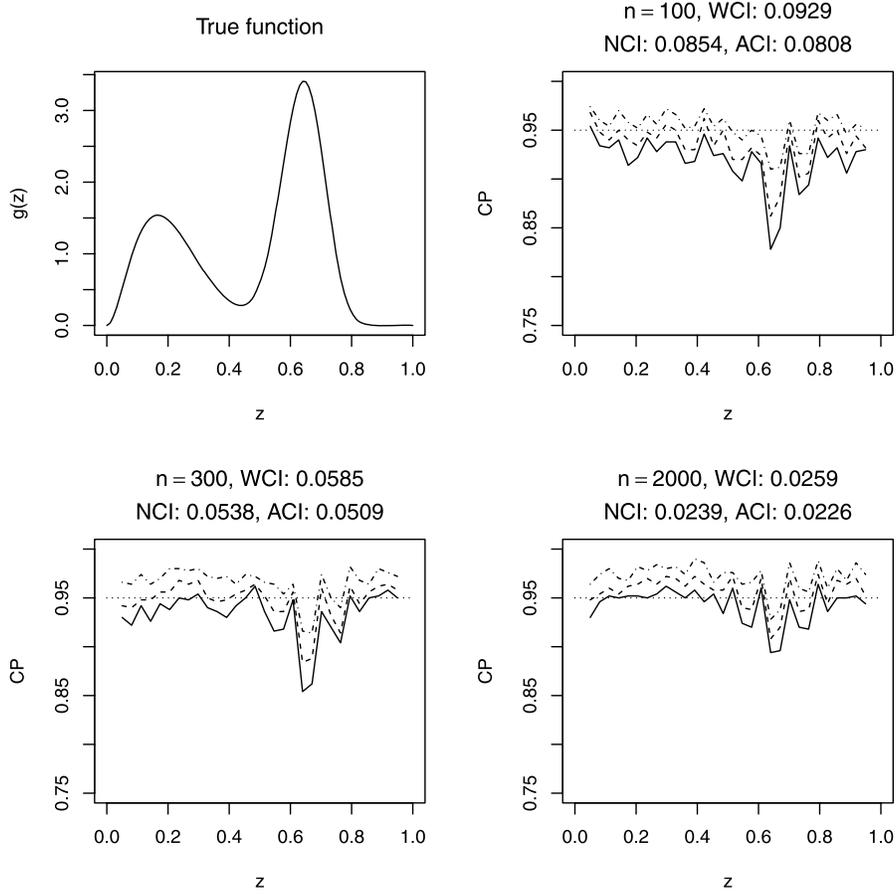}

\caption{The first panel displays the true function
$g_0(z)=0.6\beta_{30,17}(z)+0.4\beta_{3,11}(z)$ used in case~\textup{(I)} of
Example~\protect\ref{exampss}.
The other panels contain the coverage probabilities (CPs) of ACI
(solid), NCI (dashed) and WCI (dotted dashed), and the average lengths
of the three CIs (numbers in the plot titles). The CIs were built upon
thirty equally spaced covariates.}\label{CP:nonper}
\end{figure}

To examine the pointwise asymptotic CI, we considered the true function
$g_0(z)=0.6\beta_{30,17}(z)+0.4\beta_{3,11}(z)$, where $\beta_{a,b}$ is
the density function for $\operatorname{Beta}(a,b)$, and estimated it using periodic
splines with $m=2$; $\sigma$ was chosen as 0.05. In
Figure~\ref{CP:nonper}, we compare the coverage probability (CP) of our
asymptotic CI (\ref{confintformula}), denoted ACI, Wahba's Bayesian CI
(\ref{asympwahba:bci}), denoted WCI and Nychka's Bayesian CI
(\ref{asympnychka:bci}), denoted NCI, at thirty equally spaced grid
points of $\mathbb I$. The CP was computed as the proportion of the CIs
that cover $g_0$ at each point based on 1000 replications. We
observe that, in general, all CIs exhibit similar patterns, for
example, undercoverage near peaks or troughs. However, when the sample
size is sufficiently large, for example, $n=2000$, the CP of ACI is
uniformly closer to 95\% than that of WCI and NCI in smooth regions
such as $[0.1,0.4]$ and $[0.8,0.9]$. We also report the average lengths
of the three CIs in the titles of the plots. The ACI is the shortest,
as indicated in Figure~\ref{CP:nonper}.

In Figure~3 of the supplementary document \cite{SC12}, we construct the
SCB for $g$ based on formula (\ref{ss:confband}) by taking
$d_n=(-2\log{h})^{1/2}$. We compare it with the \textit{pointwise}
confidence bands constructed by linking the endpoints of the ACI, WCI
and NCI at each observed covariate, denoted ACB, BCB1 and BCB2,
respectively. The data were generated under the same setup as above. We
observe that the coverage properties of all bands are reasonably good,
and they become better as $n$ grows. Meanwhile, the band areas, that
is, the areas covered by the bands, shrink to zero as $n$ grows. We
also note that the ACB has the smallest band area, while the SCB has
the largest. This is because of the $d_n$ factor in the construction of
SCB; see Remark~\ref{remark:scb} for more details.

To conclude case~(I), we tested $H_0\dvtx g$ \textit{is linear} at the 95\%
significance level by the PLRT and GLRT. By Lemma~\ref{lemeg1} and
(\ref{H0:basis}), direct calculation leads to $r_K=1.3333$ and
$u_n=0.4714(h\sigma^{1/2})^{-1}$ when $m=2$. The data were generated
under the same setup except that different test functions
$g(z)=-0.5+z+c(\sin(\pi z)-0.5)$, $c=0,0.5,1.5,2$, were used for the
purpose of the power comparison. For the GLRT method, the R function
\textit{glkerns}(\,) provided in the \textit{lokern} package (see
\cite{E97}) was used for the local polynomial fitting based on the
Epanechnikov kernel. For the PLRT method, GCV was used to select the
smoothing parameter. Table~\ref{pow:2test} compares the power (the
proportion of rejections based on 1000 replications) for $n=20, 30,
70, 200$. When $c\geq1.5$ ($c=0$) and $n=70$ or larger, both testing
methods achieve 100\% power (5\% correct level). We also observe that
(i) the power increases as $c$ increases, that is, the test function
becomes more nonlinear; and (ii) the PLRT shows moderate advantages
over the GLRT, especially in small samples such as $n=20$. An intuitive
reason for (ii) is that the smoothing spline estimate in the PLRT uses
the full data information, while the local polynomial estimate employed
in the GLRT uses only local data information. Of course, as $n$ grows,
this difference rapidly vanishes because of the increasing data
information.
%
%t1 #&#
\begin{table}
\tabcolsep=0pt
\caption{Power comparison of the PLRT and the GLRT in case~\textup{(I)} of
Example \protect\ref{exampss} where the test function is
$g_0(z)=-0.5+z+c(\sin(\pi z)-0.5)$ with various $c$ values. The
significance level is 95\%} \label{pow:2test}
\begin{tabular*}{\tablewidth}{@{\extracolsep{\fill}}ld{2.2}d{2.2}d{2.2}d{2.2}d{3.2}d{3.2}d{3.2}d{3.2}@{}}
\hline
& \multicolumn{2}{c}{$\bolds{c=0}$} & \multicolumn{2}{c}{$\bolds{c=0.5}$}
&\multicolumn{2}{c}{$\bolds{c=1.5}$} &\multicolumn{2}{c@{}}{$\bolds{c=2}$}
\\[-4pt]
& \multicolumn{2}{c}{\hrulefill} & \multicolumn{2}{c}{\hrulefill}
&\multicolumn{2}{c}{\hrulefill} &\multicolumn{2}{c@{}}{\hrulefill}
\\
$\bolds{n}$& \multicolumn{1}{c}{\textbf{PLRT}} & \multicolumn{1}{c}{\textbf{GLRT}} &
\multicolumn{1}{c}{\textbf{PLRT}} & \multicolumn{1}{c}{\textbf{GLRT}} &
\multicolumn{1}{c}{\textbf{PLRT}} & \multicolumn{1}{c}{\textbf{GLRT}} &
\multicolumn{1}{c}{\textbf{PLRT}} & \multicolumn{1}{c@{}}{\textbf{GLRT}}
\\
\hline
& \multicolumn{8}{c@{}}{$100\times{}$Power\%}
\\
\phantom{0}20& 18.60&20.10& 28.40&30.10&89.60&86.30&97.30&96.10\\
\phantom{0}30& 13.60&14.40& 33.00&30.60&98.10&96.80&99.60&99.60\\
\phantom{0}70& 8.30 &9.40 & 54.40&48.40&100 &100 &100 &100 \\
200& 5.20 &5.50 & 95.10&92.70&100 &100 &100 &100 \\ \hline
\end{tabular*}
\end{table}

\begin{longlist}[Case (II)]
\item[Case (II). $\mathcal{H}=H^m(\mathbb{I})$:] We used
    cubic splines and repeated most of the procedures in case~(I). A
    different true function $g_0(z)=\sin(2.8\pi z)$ was chosen to
    examine the CIs. Figure~4 in the supplementary document \cite{SC12}
    summarizes the SCB and the pointwise bands constructed by ACB,
    BCB1 and BCB2. In particular, BCB1 was computed by
    (\ref{wahba:bci}) and BCB2 was constructed by scaling the length of
    BCB1 by the factor $\sqrt{27/32}\approx0.919$. We also tested the
    linearity of $g_0$ at the $95\%$ significance level, using the test
    functions $g_0(z)=-0.5+z+c(\sin(2.8\pi z)-0.5)$, for
    $c=0,0.5,1.5,2$. Table~\ref{pow:2testnp} compares the power of the
    PLRT and GLRT. From Figure~4 and Table~\ref{pow:2testnp}, we
    conclude that all findings in case~(I) are also true in case~(II).
\end{longlist}
\end{Example}

%
%ex6.2 #&#
\begin{Example}[(Nonparametric gamma model)]
Consider a two-parameter exponential model
\[
Y\mid  Z\sim\operatorname{Gamma}\bigl(\alpha,\exp\bigl(g_0(Z)\bigr)\bigr),
\]
where $\alpha>0$, $g_0\in H_0^m(\mathbb{I})$ and $Z$ is uniform over
$[0,1]$. This framework leads to $\ell(y;g(z))= \alpha
g(z)+(\alpha-1)\log{y}-y\exp(g(z))$. Thus, $I(z)=\alpha$, leading us to
choose the trigonometric polynomial basis defined as in
(\ref{H0:basis}) with $\sigma$ replaced with $\alpha ^{-1/2}$, and the
eigenvalues $\gamma_0=0$ and $\gamma_{2k}=\gamma_{2k-1}=\alpha
^{-1}(2\pi k)^{2m}$ for $k\geq1$. Local and global inference can be
conducted similarly to Example~\ref{exampss}.
\end{Example}

%
%ex6.3 #&#
\begin{Example}[(Nonparametric logistic regression)]\label{eglog}
In this example, we consider the binary response $Y\in\{0,1\}$ modeled
by the logistic relationship
%
%e6.5 #&#
\begin{equation}
\label{lgm} P(Y=1\mid Z=z)=\frac{\exp(g_0(z))}{1+\exp(g_0(z))},
\end{equation}
where $g_0\in H^m(\mathbb{I})$. Given the length of this paper, we
conducted simulations only for the ACI and PLRT. A straightforward
calculation gives $I(z)=\frac{\exp(g_0(z))}{(1+\exp(g_0(z)))^2}$, which
can be estimated by $\widehat{I}(z)=\frac{\exp(\widehat{g}_{n,\lambda
}(z))}{(1+\exp(\widehat{g}_{n,\lambda}(z)))^2}$. Given the estimate
$\widehat{I}(z)$ and the marginal density estimate $\widehat{\pi}(z)$,
we find the approximate eigenvalues and eigenfunctions via
(\ref{eigen:problem}).
%
%t2 #&#
\begin{table}%[b]
\tabcolsep=0pt
\caption{Power comparison of the PLRT and the GLRT in case~\textup{(II)} of
Example \protect\ref{exampss} where the test function is
$g_0(z)=-0.5+z+c(\sin(2.8\pi z)-0.5)$ with various $c$ values. The
significance level is 95\%} \label{pow:2testnp}
\begin{tabular*}{\tablewidth}{@{\extracolsep{\fill}}ld{2.2}d{2.2}d{3.2}d{3.2}cccc@{}}
\hline
& \multicolumn{2}{c}{$\bolds{c=0}$} & \multicolumn{2}{c}{$\bolds{c=0.5}$}
&\multicolumn{2}{c}{$\bolds{c=1.5}$} &\multicolumn{2}{c@{}}{$\bolds{c=2}$}
\\[-4pt]
& \multicolumn{2}{c}{\hrulefill} & \multicolumn{2}{c}{\hrulefill}
&\multicolumn{2}{c}{\hrulefill} &\multicolumn{2}{c@{}}{\hrulefill}
\\
$\bolds{n}$& \multicolumn{1}{c}{\textbf{PLRT}} & \multicolumn{1}{c}{\textbf{GLRT}} & \multicolumn{1}{c}{\textbf{PLRT}}
   & \multicolumn{1}{c}{\textbf{GLRT}} & \multicolumn{1}{c}{\textbf{PLRT}} & \textbf{GLRT} & \textbf{PLRT} & \textbf{GLRT}
   \\
\hline
& \multicolumn{8}{c@{}}{$100\times{}$Power\%}
\\
\phantom{0}20& 16.00&17.40&71.10&67.60&100&100&100&100\\
\phantom{0}30& 12.70&14.00&83.20&81.20&100&100&100&100\\
\phantom{0}70& 6.50 &7.40 &99.80&99.70&100&100&100&100\\
200& 5.10 &5.30 &100 &100 &100&100&100&100\\ \hline
\end{tabular*}
\end{table}

The results are based on 1000 replicated data sets drawn from
(\ref{lgm}), with $n=70,100,300,500$. To test whether $g$ is linear, we
considered two test functions, $g_0(z)=-0.5+z+c(\sin(\pi z)-0.5)$ and
$g_0(z)=-0.5+z+c(\sin(2.8\pi z)-0.5)$, for $c=0,1,1.5,2$. We use $m=2$.
Numerical calculations reveal that the eigenvalues are
$\gamma_\nu\approx(\alpha\nu)^{2m}$, where $\alpha
=4.40,4.41,4.47,4.52$ and $\alpha=4.40,4.44,4.71,4.91$ corresponding to
the two test functions and the four values of $c$. This simplifies the
calculations of $\sigma_K^2$ and $\rho_K^2$ defined in Theorem
\ref{Global:LRT:Lim}. For instance, when
$\gamma_\nu\approx(4.40\nu)^{2m}$, using a result analogous to Lemma
\ref{lemeg1} we have $\sigma_K^2\approx0.25$ and $\rho_K^2\approx0.19$.
Then the quantities $r_K$ and $u_n$ are found for the PLRT method. To
evaluate ACI, we considered the true function
$g_0(z)=(0.15)10^6z^{11}(1-z)^6+(0.5)10^4z^3(1-z)^{10}-1$. The CP and
the average lengths of the ACI are calculated at thirty evenly spaced
points in $\mathbb I$ under three sample sizes, $n=200, 500, 2000$.
%
%t3 #&#
\begin{table}[t]
\tabcolsep=0pt
\tablewidth=220pt
\caption{Power of PLRT in Example~\protect\ref{eglog} where the test
function is $g_0(z)=-0.5+z+c(\sin(\pi z)-0.5)$ with various $c$ values.
The significance level is 95\%} \label{pow:1glm}
\begin{tabular*}{\tablewidth}{@{\extracolsep{\fill}}lcccd{3.2}@{}}
\hline
$\bolds{n}$& $\bolds{c=0}$ & $\bolds{c=1}$ & $\bolds{c=1.5}$ & \multicolumn{1}{c@{}}{$\bolds{c=2}$}\\
\hline
& \multicolumn{4}{c@{}}{$100\times{}$Power\%}
\\
\phantom{0}70 & 4.10 & 16.90 & 30.20 &50.80\\
100& 4.50 & 17.30 & 38.90 &63.40\\
300& 5.00 & 52.50 & 92.00 &99.30\\
500& 5.00 & 79.70 & 99.30 &100 \\
\hline
\end{tabular*}
\end{table}
%
%t4 #&#
\begin{table}
\tablewidth=220pt
\tabcolsep=0pt
\caption{Power of PLRT in Example~\protect\ref{eglog} where the test
function is $g_0(z)=-0.5+z+c(\sin(2.8\pi z)-0.5)$ with various $c$
values. The significance level is 95\%} \label{pow:2glm}
\begin{tabular*}{\tablewidth}{@{\extracolsep{\fill}}lcd{3.2}d{3.2}d{3d2}@{}}
\hline
$\bolds{n}$& $\bolds{c=0}$ & \multicolumn{1}{c}{$\bolds{c=1}$} & \multicolumn{1}{c}{$\bolds{c=1.5}$} & \multicolumn{1}{c@{}}{$\bolds{c=2}$}
\\
\hline
& \multicolumn{4}{c@{}}{$100\times{}$Power\%}
\\
\phantom{0}70 &4.10 & 56.20 & 90.10 &99.00\\
100&5.00 & 71.90 & 96.90 &100\\
300&5.00 & 99.80 & 100 &100\\
500&5.00 & 100 & 100 &100\\
\hline
\end{tabular*}
\end{table}

%
%f2 #&#
\begin{figure}%[t]

\includegraphics{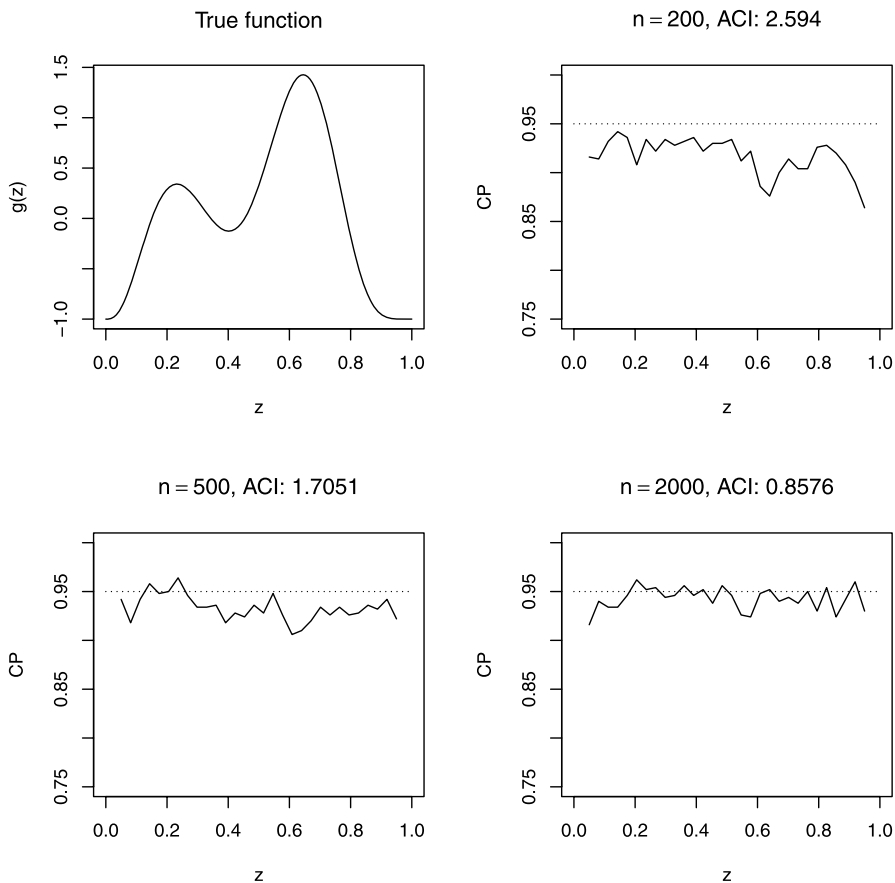}

\caption{The first panel displays the true function
$g_0(z)=(0.15)10^6z^{11}(1-z)^6+\allowbreak(0.5)10^4z^3(1-z)^{10}-1$ used in
Example~\protect\ref{eglog}.
The other panels contain the CP and average length (number in the plot title)
of each ACI. The ACIs were built upon thirty equally spaced
covariates.}\label{CP:logit}\vspace*{6pt}
\end{figure}

The results on the power of the PLRT are summarized in
Tables~\ref{pow:1glm} and~\ref{pow:2glm}, which demonstrate the
validity of the proposed testing method. Specifically, when \mbox{$c=0$}, the
power reduces to the desired size 0.05; when $c\geq1.5$ and
$n\geq300$, the power approaches one. The results for the CPs and
average lengths of ACIs are summarized in Figure~\ref{CP:logit}. The CP
uniformly approaches the desired 95\% confidence level as $n$ grows,
showing the validity of the intervals.
\end{Example}

\section*{Acknowledgement}
We are grateful for helpful discussions with
Professor Chong Gu and careful proofreading of Professor Pang Du. The
authors also thank the Co-Editor Peter Hall, the Associate Editor, and
two referees for insightful comments that led to important improvements
in the paper.

\begin{supplement}%[id=suppA]
\stitle{Supplement to ``Local and global asymptotic inference in~smoothing spline models''}
\slink[doi]{10.1214/13-AOS1164SUPP} %[doi,text={...}] - jei reikia
%suskaldyti doi
\sdatatype{.pdf}
\sfilename{aos1164\_supp.pdf}
\sdescription{The supplementary materials contain all the proofs of the theoretical results in the present paper.}
\end{supplement}

% imsref loaded by linak, 2013-09-27 16:29:12
% imsref loaded by linak, 2013-10-01 14:47:06

\printaddresses

\end{document}